\documentclass[12pt]{amsart}

\usepackage[T1]{fontenc}
\usepackage[utf8]{inputenc}
\usepackage{microtype}
\usepackage{cite}

\usepackage{amsmath}
\usepackage{amssymb}
\usepackage{mathtools}
\usepackage{mathrsfs}
\usepackage{amsthm}
\usepackage{fullpage}

\usepackage{enumerate}
\usepackage{tikz}
\usepackage{xcolor}

\usepackage[
    colorlinks=true,
    linkcolor=blue,
    citecolor=blue,
    urlcolor=blue
]{hyperref}

\usepackage[
    capitalize,
    nameinlink
]{cleveref}

\theoremstyle{plain}

\newtheorem{theorem}{Theorem}[section]
\newtheorem{lemma}[theorem]{Lemma}
\newtheorem{proposition}[theorem]{Proposition}
\newtheorem{corollary}[theorem]{Corollary}

\newtheorem*{theorem*}{Theorem}

\theoremstyle{definition}

\theoremstyle{remark}

\theoremstyle{plain}

\newtheorem{manualtheoreminner}{Theorem}

\numberwithin{equation}{section}

\title[Sharp constants and asymptotics for Riesz potentials]
{Sharp Constants and First-Order Asymptotics in a
Bourgain--Brezis--Mironescu Formula for Riesz Potentials}

\author[M. Jleli]{Mohamed Jleli}

\address{(M. Jleli) Department of Mathematics,
College of Science,
King Saud University,
Riyadh 11451, Saudi Arabia}

\email{jleli@ksu.edu.sa}

\author[B. Samet]{Bessem Samet}

\address{(B. Samet) Department of Mathematics,
College of Science,
King Saud University,
Riyadh 11451, Saudi Arabia}

\email{bessem.samet@gmail.com}

\date{}



\begin{document}


\begin{abstract}
We investigate sharp constants and asymptotic properties of a pointwise inequality of Bourgain--Brezis--Mironescu type in dimension \(n\geq2\), involving Riesz potentials and Spector's nonlinear fractional differential operator \(\mathcal D^\alpha\) of order \(0<\alpha<1\). For every \(n\geq3\) and \(0<\alpha<1\), we determine the optimal constant without symmetry restrictions and prove that it is attained by smooth radially monotone functions centered at the evaluation point. The proof combines geodesic-flow invariance, kernel monotonicity, and circular rearrangement. In dimension two, we establish the leading-order asymptotic behavior of the optimal constant and derive an explicit formula in the radial class. We further refine the pointwise limit established by A. Claros and C. P\'erez [Bourgain--Brezis--Mironescu formula for Riesz potentials, Adv. Math. 503 (2026), Art. 111226] by obtaining an explicit first-order correction with a quadratic remainder uniform in the spatial variable. 
Finally, we strengthen their convergence-in-measure result on compact subsets for functions in
\(W^{1,1}(\mathbb R^n)\) by establishing almost-everywhere convergence as \(\alpha\to1^-\), without passing to a subsequence.

\end{abstract}

\maketitle


\section{Introduction and main results}

In their seminal works \cite{BBM1,BBM2}, Bourgain, Brezis, and Mironescu established that suitably normalized fractional Gagliardo energies converge to their first-order Sobolev counterparts as the fractional order approaches one from below. Since then, BBM-type limiting phenomena have been investigated in various settings, including fractional Orlicz--Sobolev spaces \cite{AlbericoCianchiPickSlavikova}, interpolation and extrapolation frameworks \cite{DominguezMilman}, whole-space convergence results involving Triebel--Lizorkin spaces \cite{BrazkeSchikorraYung}, generalized fractional Sobolev seminorms on arbitrary domains \cite{Mohanta}, and ball Banach function spaces \cite{DaiEtAl}. Related characterizations of ball Banach Sobolev spaces can be found in \cite{ZhuYangYuan}.

In particular, for \(n\geq1\) and
\(f\in C_c^\infty(\mathbb R^n)\), the BBM formula
in the case \(p=1\), which is central to the
present work, reads
\[
\lim_{\alpha\to1^-}
(1-\alpha)
\int_{\mathbb R^n}\int_{\mathbb R^n}
\frac{|f(x)-f(y)|}{|x-y|^{n+\alpha}}
\,dx\,dy
=
K_n
\int_{\mathbb R^n}|\nabla f(x)|\,dx,
\]
where
\[
K_n
=
\int_{\mathbb S^{n-1}}
|\omega\cdot e|\,d\sigma(\omega)
=
\frac{\Gamma(n/2)}
{\sqrt{\pi}\,\Gamma((n+1)/2)}
|\mathbb S^{n-1}|,
\qquad e\in\mathbb S^{n-1},
\]
\(d\sigma\) denotes the standard surface measure
on \(\mathbb S^{n-1}\), and \(\Gamma\) is the
Euler Gamma function. The factor \(1-\alpha\)
compensates for the divergence of the fractional
energy as \(\alpha\to1^-\).

Quantitative and weighted refinements of BBM-type inequalities have also been investigated. Weighted fractional Poincar\'e--Sobolev and Hardy-type inequalities, with explicit control of the constants as the fractional parameter approaches one, were established in \cite{HurriSyrjanenEtAl}. More recently, Hu, Li, Yang, and Yuan \cite{HuLiYangYuan} established localized weighted inequalities relating Gagliardo and Sobolev seminorms and derived BBM-type formulas in the framework of ball Banach function spaces.

Convergence rates and higher-order asymptotics in BBM-type formulas have also been investigated. Chambolle, Novaga, and Pagliari \cite{ChambolleNovagaPagliari} studied convergence rates for nonlocal energies and established the \(\Gamma\)-convergence of suitably rescaled functionals to a second-order limit functional. Second-order asymptotic expansions were obtained for fractional perimeters in \cite{CesaroniNovaga} and, more recently, for quadratic fractional Gagliardo energies in \cite{KubinPagliariTribuzio}. Fern\'andez Bonder and Salort \cite{FernandezBonderSalort} investigated anisotropic fractional energies, showing how the directional dependence of nonlocal interactions is reflected in the limiting local energy.

Riesz potentials provide a natural framework for Sobolev-type subrepresentation formulas. For \(0<\beta<n\) and a suitable function \(g\), let
\[
I_\beta g(x)
=
\gamma_{n,\beta}
\int_{\mathbb R^n}
\frac{g(y)}{|x-y|^{n-\beta}}\,dy,
\]
where
\begin{equation}\label{eq:gamma-n-beta}
\gamma_{n,\beta}
=
\frac{\Gamma\bigl((n-\beta)/2\bigr)}
{2^\beta\pi^{n/2}\Gamma(\beta/2)}.
\end{equation}
For \(n\geq2\), the classical subrepresentation estimate
\[
|f(x)|
\lesssim_n
I_1(|\nabla f|)(x),
\qquad
f\in C_c^\infty(\mathbb R^n),
\]
provides pointwise control of a function by a Riesz potential of its gradient. Related representation formulas and weighted Poincar\'e inequalities have been studied in more general settings; see, for instance, \cite{FranchiLuWheeden}. Haj\l asz and Liu \cite{HajlaszLiu} investigated maximal potentials and their Sobolev regularity. More recently, Hoang, Moen, and P\'erez established pointwise estimates for rough operators in terms of Riesz potentials \cite{HoangMoenPerezRough}, as well as related estimates involving sparse potential operators \cite{HoangMoenPerezJFA}.

A fractional analogue of this subrepresentation
framework involves the nonlinear fractional
differential operator introduced by Spector
\cite{Spector},
\[
\mathcal D^\alpha f(x)
=
\int_{\mathbb R^n}
\frac{|f(x)-f(y)|}{|x-y|^{n+\alpha}}\,dy,
\qquad
0<\alpha<1.
\]
In connection with the BBM formula, it is natural
to consider the normalized Riesz potential
\[
(1-\alpha)I_\alpha(\mathcal D^\alpha f)
\]
and investigate its relationship with
\(I_1(|\nabla f|)\) as \(\alpha\to1^-\).

Hoang, Moen, and P\'erez
\cite{HoangMoenPerezSubrepresentation}
established that, for \(n\geq2\), \(0<\alpha<1\),
and \(f\in C_c^\infty(\mathbb R^n)\),
\begin{equation}\label{eq:HMP-comparison}
(1-\alpha)I_\alpha(\mathcal D^\alpha f)(x)
\leq
C_{n,\alpha}^{\mathrm{HMP}} I_1(|\nabla f|)(x),
\qquad x\in\mathbb R^n,
\end{equation}
where, with the normalization
\eqref{eq:gamma-n-beta},
\begin{equation}\label{eq-C-HMP}
C_{n,\alpha}^{\mathrm{HMP}}
=
\frac{
2^{2-\alpha}\pi^{n/2}
\Gamma\bigl((3-\alpha)/2\bigr)
}{
\alpha\,
\Gamma\bigl((n+\alpha-1)/2\bigr)
}.
\end{equation}
In particular,
\[
\lim_{\alpha\to1^-}C_{n,\alpha}^{\mathrm{HMP}}
=
|\mathbb S^{n-1}|.
\]
The constant \(C_{n,\alpha}^{\mathrm{HMP}}\)
is explicit, but its optimality was not
addressed in
\cite{HoangMoenPerezSubrepresentation}.

Claros and P\'erez \cite{ClarosPerez} established
the pointwise BBM limit of the normalized Riesz
potential. More precisely, for \(n\geq2\),
\(f\in C_c^\infty(\mathbb R^n)\), and every
\(x\in\mathbb R^n\),
\begin{equation}\label{eq:Claros-Perez-limit}
\lim_{\alpha\to1^-}
(1-\alpha)I_\alpha(\mathcal D^\alpha f)(x)
=
K_n I_1(|\nabla f|)(x).
\end{equation}
They further proved, by a density argument,
that for \(f\in W^{1,1}(\mathbb R^n)\), the same
convergence holds in measure on every compact
subset of \(\mathbb R^n\). Consequently, every
sequence \(\alpha_k\to1^-\) admits a subsequence
along which the convergence holds almost
everywhere in \(\mathbb R^n\).

For \(n\geq2\), \(0<\alpha<1\),
\(f\in C_c^\infty(\mathbb R^n)\setminus\{0\}\),
and \(x\in\mathbb R^n\), we consider the quotient
\begin{equation}\label{eq:R-alpha}
R_\alpha(f,x)
=
\frac{
(1-\alpha)I_\alpha(\mathcal D^\alpha f)(x)
}{
I_1(|\nabla f|)(x)
}.
\end{equation}
We define the optimal constant associated with
\eqref{eq:R-alpha} by
\begin{equation}\label{eq:def-C-star}
C_{n,\alpha}^{*}
=
\sup_{\substack{
f\in C_c^\infty(\mathbb R^n)\setminus\{0\}\\
x\in\mathbb R^n}}
R_\alpha(f,x).
\end{equation}
Equivalently, \(C_{n,\alpha}^{*}\) is the smallest
constant such that
\[
(1-\alpha)I_\alpha(\mathcal D^\alpha f)(x)
\leq
C_{n,\alpha}^{*} I_1(|\nabla f|)(x)
\]
for every \(f\in C_c^\infty(\mathbb R^n)\)
and every \(x\in\mathbb R^n\).
By \eqref{eq:HMP-comparison},
\[
0<C_{n,\alpha}^{*}
\leq
C_{n,\alpha}^{\mathrm{HMP}}
<\infty.
\]

For each fixed admissible pair \((f,x)\),
the pointwise limit \eqref{eq:Claros-Perez-limit}
yields
\[
\lim_{\alpha\to1^-}R_\alpha(f,x)=K_n.
\]
However, this convergence does not justify
interchanging the limit with the supremum
in \eqref{eq:def-C-star}, since maximizing
pairs may depend on \(\alpha\).

The first objective of the present paper is to study the optimal constant \(C_{n,\alpha}^{*}\). In dimension two, we establish its leading-order asymptotic behavior as \(\alpha\to1^-\), determine the exact constant in the radial class for every \(0<\alpha<1\), and identify a class of radial extremizers. In dimensions \(n\geq3\), we determine the exact value of \(C_{n,\alpha}^{*}\) for every \(0<\alpha<1\) and identify a class of extremizers. The second objective is to refine the leading-order BBM limit \eqref{eq:Claros-Perez-limit} by deriving a first-order asymptotic expansion of the normalized potential \((1-\alpha)I_\alpha(\mathcal D^\alpha f)\), with a remainder estimate uniform in \(x\in\mathbb R^n\). Finally, we strengthen the subsequential almost-everywhere convergence result of Claros and P\'erez \cite{ClarosPerez} by proving convergence almost everywhere as \(\alpha\to1^-\), without passing to a subsequence, for \(f\in W^{1,1}(\mathbb R^n)\).

We begin with the planar case. Our first result
describes the leading-order asymptotic behavior
of the unrestricted optimal constant
\(C_{2,\alpha}^{*}\).

\begin{theorem}\label{thm:leading-asymptotics-2d}
As \(\alpha\to1^-\),
\begin{equation}\label{eq:leading-asymptotics-C2}
C_{2,\alpha}^{*}
=
4+O(1-\alpha).
\end{equation}
\end{theorem}

In particular, \eqref{eq:leading-asymptotics-C2}
yields
\[
\lim_{\alpha\to1^-}C_{2,\alpha}^{*}=4=K_2,
\]
in agreement with the pointwise BBM limit
\eqref{eq:Claros-Perez-limit}.

We next restrict the optimization problem to radially symmetric functions centered at the evaluation point. For \(0<\alpha<1\), define
\[
C_{2,\alpha}^{\mathrm{rad}}
=
\sup_{\substack{
x_0\in\mathbb R^2\\
f(y)=F(|y-x_0|)\in C_c^\infty(\mathbb R^2)\setminus\{0\}
}}
R_\alpha(f,x_0).
\]
Our next result provides an explicit formula for this constant and identifies a class of extremizers.

\begin{theorem}\label{thm:sharp-constant-radial-2d}
For every \(0<\alpha<1\),
\begin{equation}\label{eq:sharp-C2-rad}
C_{2,\alpha}^{\mathrm{rad}}
=
\frac{2^{2-\alpha}\pi(1-\alpha)}{\alpha^2}
\frac{\Gamma(1-\alpha/2)}{\Gamma(\alpha/2)}
\left(
2\frac{\Gamma(1-\alpha)}
{\Gamma(1-\alpha/2)^2}-1
\right).
\end{equation}
Moreover, for every \(x_0\in\mathbb R^2\) and every nontrivial function
\[
f(y)=F(|y-x_0|)\in C_c^\infty(\mathbb R^2),
\]
with \(F\) monotone, we have
\[
R_\alpha(f,x_0)=C_{2,\alpha}^{\mathrm{rad}}.
\]
\end{theorem}

Theorem~\ref{thm:sharp-constant-radial-2d} provides an explicit formula for the sharp constant in the radial class and identifies a family of extremizers for every \(0<\alpha<1\). The determination of the unrestricted optimal constant \(C_{2,\alpha}^{*}\) for arbitrary \(0<\alpha<1\) remains unresolved by the present analysis.

We now turn to the case \(n\geq3\). For \(0<\alpha<1\), define
\begin{equation}\label{eq:def-radial-coefficient-n}
\mathfrak C_{n,\alpha}^{\rm rad}
=
|\mathbb S^{n-2}|
\int_0^\infty \tau
\int_0^\pi
\frac{
2\cosh\!\left(\frac{n-\alpha}{2}\tau\right)
(\sin\theta)^{n-2}
}{
\bigl[2(\cosh\tau-\cos\theta)\bigr]^{(n+\alpha)/2}
}
\,d\theta\,d\tau.
\end{equation}
The integral in \eqref{eq:def-radial-coefficient-n}
is finite for every \(0<\alpha<1\).

Our next main result is the following.

\begin{theorem}\label{thm:sharp-constant-higher-dim}
Let \(n\geq3\). For every \(0<\alpha<1\),
\begin{equation}\label{eq:sharp-Cn}
C_{n,\alpha}^{*}
=
(1-\alpha)
\frac{\gamma_{n,\alpha}}{\gamma_{n,1}}
\mathfrak C_{n,\alpha}^{\rm rad}.
\end{equation}
Moreover, for every \(x_0\in\mathbb R^n\) and every nontrivial function
\[
f(y)=F(|y-x_0|)\in C_c^\infty(\mathbb R^n),
\]
with \(F\) monotone, we have
\[
R_\alpha(f,x_0)=C_{n,\alpha}^{*}.
\]
\end{theorem}

Theorem~\ref{thm:sharp-constant-higher-dim}
shows that, for \(n\geq3\), the unrestricted
optimization problem admits an exact reduction
to the class of radially monotone functions.
In fact, any fixed nontrivial radially
monotone function in \(C_c^\infty(\mathbb R^n)\),
centered at the origin, attains the optimal
constant for every \(0<\alpha<1\).
Consequently, by \eqref{eq:Claros-Perez-limit},
\[
\lim_{\alpha\to1^-}C_{n,\alpha}^{*}=K_n.
\]

The following corollary shows that the sharp
constant provides a strict improvement over
the constant in the Hoang--Moen--P\'erez
inequality \eqref{eq:HMP-comparison} for every
\(0<\alpha<1\), and quantifies this improvement
as \(\alpha\to1^-\).

\begin{corollary}\label{cor:strict-HMP}
Let \(n\geq3\). Then, for every \(0<\alpha<1\),
\[
C_{n,\alpha}^{*}
<
C_{n,\alpha}^{\mathrm{HMP}}.
\]
Moreover,
\[
\lim_{\alpha\to1^-}
\frac{C_{n,\alpha}^{*}}
{C_{n,\alpha}^{\mathrm{HMP}}}
=
\frac{K_n}{|\mathbb S^{n-1}|}<1.
\]
\end{corollary}

Our next main result refines the pointwise BBM
limit \eqref{eq:Claros-Perez-limit} by identifying
the first-order asymptotic correction for each fixed
\(f\in C_c^2(\mathbb R^n)\) and establishing a
remainder estimate uniform in \(x\in\mathbb R^n\).

For \(n\geq2\), \(f\in C_c^2(\mathbb R^n)\), and
\(y\in\mathbb R^n\), define
\begin{equation}\label{eq:def-Ff}
\mathcal F_f(y)
=
\int_{|h|<1}
\frac{
|f(y+h)-f(y)|
-
|\nabla f(y)\cdot h|
}{
|h|^{n+1}
}\,dh
+
\int_{|h|\geq1}
\frac{|f(y+h)-f(y)|}{|h|^{n+1}}\,dh.
\end{equation}
Both integrals in \eqref{eq:def-Ff} are absolutely
convergent for every \(y\in\mathbb R^n\).
Moreover,
\[
\mathcal F_f\in L^1(\mathbb R^n)\cap L^\infty(\mathbb R^n).
\]

Let
\begin{equation}\label{eq:def-lambda-n}
\lambda_n
=
-\frac12\psi\!\left(\frac{n-1}{2}\right)
-\log 2
-\frac12\psi\!\left(\frac12\right),
\end{equation}
where \(\psi=\Gamma'/\Gamma\) denotes the digamma function.

For \(x\in\mathbb R^n\), we then define
\begin{equation}\label{eq:def-Af}
\mathcal A_f(x)
=
I_1(\mathcal F_f)(x)
-
K_n\gamma_{n,1}
\int_{\mathbb R^n}
\frac{\log|x-y|+\lambda_n}{|x-y|^{n-1}}
|\nabla f(y)|\,dy.
\end{equation}
Both terms on the right-hand side of
\eqref{eq:def-Af} are absolutely convergent
for every \(x\in\mathbb R^n\). Moreover,
\[
\mathcal A_f\in L^\infty(\mathbb R^n).
\]

\begin{theorem}\label{thm:first-order-expansion}
Let \(n\geq2\) and \(f\in C_c^2(\mathbb R^n)\). Then, as \(\alpha\to1^-\),
\begin{equation}\label{eq:first-order-expansion}
\left\|
(1-\alpha)I_\alpha(\mathcal D^\alpha f)
-
K_n I_1(|\nabla f|)
-
(1-\alpha)\mathcal A_f
\right\|_{L^\infty(\mathbb R^n)}
=
O\bigl((1-\alpha)^2\bigr).
\end{equation}
\end{theorem}

The correction term \(\mathcal A_f\) consists of two contributions: \(I_1(\mathcal F_f)\) arises from the finite remainder in the expansion of the nonlinear fractional differential operator \(\mathcal D^\alpha f\), while the logarithmic term in \eqref{eq:def-Af} reflects the first-order variation of the Riesz kernel and its normalization as \(\alpha\to1^-\).

Our final main result strengthens the
\(W^{1,1}\)-extension of Claros and P\'erez
\cite{ClarosPerez} by establishing almost
everywhere convergence as \(\alpha\to1^-\),
without passing to a subsequence.

\begin{theorem}\label{thm:W11-full-limit}
Let \(n\geq2\) and \(f\in W^{1,1}(\mathbb R^n)\).
Then, for almost every \(x\in\mathbb R^n\),
\[
\lim_{\alpha\to1^-}
(1-\alpha)I_\alpha(\mathcal D^\alpha f)(x)
=
K_n I_1(|\nabla f|)(x).
\]
\end{theorem}

We now briefly outline the main ideas underlying the proofs of our results. Our approach to the sharp-constant problem is based on logarithmic-polar coordinates, which allow us to reformulate the problem as an integral inequality on the cylinder \(\mathbb R\times\mathbb S^{n-1}\).

In dimension two, the main challenge is to
control the directional dependence of the
variational upper bound as \(\alpha\to1^-\).
We overcome this difficulty by establishing
a uniform asymptotic estimate for the associated
directional functional, which yields
Theorem~\ref{thm:leading-asymptotics-2d}.
Furthermore, we solve the optimization problem
in the radial class for every \(0<\alpha<1\),
obtaining an explicit formula for the sharp
radial constant and identifying a family of
extremizers, as stated in
Theorem~\ref{thm:sharp-constant-radial-2d}.

In dimensions \(n\geq3\), the main challenge is to establish the global optimality of the radial direction for every \(0<\alpha<1\). Exploiting the invariance of the natural measure under the geodesic flow on the unit tangent bundle of \(\mathbb S^{n-1}\), we derive an upper bound for the unrestricted problem in terms of a directional maximization problem. The key step is to prove that the associated directional functional attains its maximum in the radial direction. This is achieved by establishing a monotonicity property of a suitably transformed kernel and applying a circular rearrangement inequality. The resulting sharp bound is attained by radially monotone functions, yielding Theorem~\ref{thm:sharp-constant-higher-dim}.

Turning to the first-order asymptotics, the main difficulty is to identify the correction arising from the simultaneous variation of the nonlinear fractional operator \(\mathcal D^\alpha\) and the Riesz kernel as \(\alpha\to1^-\). Our approach consists in isolating the leading singular contribution of \(\mathcal D^\alpha f\) and controlling its remainder in both \(L^1\) and \(L^\infty\). We then combine this analysis with a first-order expansion of the normalized Riesz kernel, including its logarithmic correction. This allows us to identify the explicit correction term \(\mathcal A_f\) and establish the uniform quadratic remainder estimate in Theorem~\ref{thm:first-order-expansion}.

Finally, the proof of Theorem~\ref{thm:W11-full-limit} relies on a pointwise stability estimate, uniform with respect to \(\alpha\) near \(1\), derived from the Hoang--Moen--P\'erez inequality \eqref{eq:HMP-comparison}. Combined with a suitable approximation by smooth functions in \(W^{1,1}(\mathbb R^n)\), this estimate allows us to construct a common set of full measure, independent of \(\alpha\), and thereby establish almost-everywhere convergence without passing to a subsequence.

The rest of the paper is organized as follows. Section~\ref{sec2} is devoted to the planar case, where we prove Theorems~\ref{thm:leading-asymptotics-2d} and~\ref{thm:sharp-constant-radial-2d}. Section~\ref{sec3} contains the proof of Theorem~\ref{thm:sharp-constant-higher-dim} and the derivation of Corollary~\ref{cor:strict-HMP}. Theorem~\ref{thm:first-order-expansion} is proved in Section~\ref{sec4}, while the proof of Theorem~\ref{thm:W11-full-limit} is given in Section~\ref{sec5}. Finally, Section~\ref{sec6} discusses the question of determining the unrestricted sharp constant in dimension two and establishes an explicit series representation for the sharp constant in higher dimensions.

\section{The constants \(C_{2,\alpha}^{*}\) and \(C_{2,\alpha}^{\mathrm{rad}}\)}\label{sec2}

This section is devoted to the proofs of
Theorems~\ref{thm:leading-asymptotics-2d}
and~\ref{thm:sharp-constant-radial-2d}.
We begin with several auxiliary lemmas.

For \(0<\alpha<1\) and
\[
(\tau,\vartheta)\in
\bigl(\mathbb R\times[-\pi,\pi]\bigr)\setminus\{(0,0)\},
\]
define
\begin{equation}\label{eq:def-K-alpha-2d}
\mathscr K_\alpha(\tau,\vartheta)
=
\frac{
2\cosh\!\left(\frac{2-\alpha}{2}\tau\right)
}{
\bigl[2(\cosh\tau-\cos\vartheta)\bigr]^{1+\alpha/2}
}.
\end{equation}
For \(\beta\in[0,\pi/2]\), set
\begin{equation}\label{eq:def-A-alpha-beta}
A_\alpha(\beta)
=
2\int_0^\infty\int_0^\pi
\max\bigl\{
\tau\cos\beta,\,
\vartheta\sin\beta
\bigr\}
\mathscr K_\alpha(\tau,\vartheta)
\,d\vartheta\,d\tau.
\end{equation}
The integral defining \(A_\alpha(\beta)\) is finite
for every \(0<\alpha<1\) and
\(\beta\in[0,\pi/2]\).

\begin{lemma}\label{L2.1}
For every \(0<\alpha<1\),
\begin{equation}\label{eq:C2-log-polar-upper}
C_{2,\alpha}^{*}
\leq
(1-\alpha)
\frac{\gamma_{2,\alpha}}{\gamma_{2,1}}
\sup_{0\leq\beta\leq\pi/2}A_\alpha(\beta).
\end{equation}
\end{lemma}

\begin{proof}
Fix \(0<\alpha<1\). By translation invariance,
\eqref{eq:def-C-star} gives
\[
C_{2,\alpha}^{*}
=
\sup_{f\in C_c^\infty(\mathbb R^2)\setminus\{0\}}
R_\alpha(f,0).
\]

Let \(f\in C_c^\infty(\mathbb R^2)\setminus\{0\}\). We introduce
logarithmic-polar coordinates
\[
y=e^t\omega(\theta),
\qquad
\omega(\theta)=(\cos\theta,\sin\theta),
\]
and set
\[
u(t,\theta)
=
f\bigl(e^t\omega(\theta)\bigr),
\qquad
(t,\theta)\in\mathbb R\times\mathbb S^1.
\]
Then
\[
|\nabla f(e^t\omega(\theta))|
=
e^{-t}
\sqrt{
|\partial_tu(t,\theta)|^2+
|\partial_\theta u(t,\theta)|^2
},
\]
and therefore
\begin{equation}\label{eq:denominator-cylinder-2d}
\frac{1}{\gamma_{2,1}}
I_1(|\nabla f|)(0)
=
\int_{\mathbb R}\int_0^{2\pi}
|\nabla_{\!c}u(t,\theta)|
\,d\theta\,dt,
\end{equation}
where
\[
\nabla_{\!c}u
=
(\partial_tu,\partial_\theta u)
\]
denotes the gradient on the cylinder \(\mathbb R\times\mathbb S^1\).

Writing
\[
y=e^t\omega(\theta),
\qquad
z=e^s\omega(\varphi),
\]
we have
\[
|y-z|^2
=
2e^{t+s}
\bigl(\cosh(t-s)-\cos(\theta-\varphi)\bigr).
\]
Since
\[
dy=e^{2t}\,dt\,d\theta,
\qquad
dz=e^{2s}\,ds\,d\varphi,
\]
it follows from the definitions of \(I_\alpha\) and
\(\mathcal D^\alpha\) that
\[
\begin{aligned}
\frac{1}{\gamma_{2,\alpha}}
I_\alpha(\mathcal D^\alpha f)(0)
=
\int_{\mathbb R}\int_0^{2\pi}
\int_{\mathbb R}\int_0^{2\pi}
\frac{
e^{\alpha t+2s}
|u(t,\theta)-u(s,\varphi)|
}{
\bigl[
2e^{t+s}
\bigl(\cosh(t-s)-\cos(\theta-\varphi)\bigr)
\bigr]^{1+\alpha/2}
}
\,d\varphi\,ds\,d\theta\,dt.
\end{aligned}
\]
With the change of variables
\[
\tau=s-t,
\qquad
\vartheta=\varphi-\theta,
\]
and using the \(2\pi\)-periodicity in the angular variable, we obtain
\[
\begin{aligned}
&\frac{1}{\gamma_{2,\alpha}}
I_\alpha(\mathcal D^\alpha f)(0)
\\
&\quad =
\int_{\mathbb R}\int_{-\pi}^{\pi}
\frac{
e^{(1-\alpha/2)\tau}
}{
\bigl[2(\cosh\tau-\cos\vartheta)\bigr]^{1+\alpha/2}
}
\Bigg[
\int_{\mathbb R}\int_0^{2\pi}
\left|
u(t+\tau,\theta+\vartheta)-u(t,\theta)
\right|
\,d\theta\,dt
\Bigg]
\,d\vartheta\,d\tau.
\end{aligned}
\]
On the other hand, interchanging the pairs \((t,\theta)\) and
\((s,\varphi)\) in the original integral and performing the same change
of variables gives
\[
\begin{aligned}
&\frac{1}{\gamma_{2,\alpha}}
I_\alpha(\mathcal D^\alpha f)(0)
\\
&\quad =
\int_{\mathbb R}\int_{-\pi}^{\pi}
\frac{
e^{-(1-\alpha/2)\tau}
}{
\bigl[2(\cosh\tau-\cos\vartheta)\bigr]^{1+\alpha/2}
}
\Bigg[
\int_{\mathbb R}\int_0^{2\pi}
\left|
u(t+\tau,\theta+\vartheta)-u(t,\theta)
\right|
\,d\theta\,dt
\Bigg]
\,d\vartheta\,d\tau.
\end{aligned}
\]
Averaging the two identities and using
\eqref{eq:def-K-alpha-2d}, we obtain
\begin{equation}\label{eq:numerator-cylinder-2d}
\frac{1}{\gamma_{2,\alpha}}
I_\alpha(\mathcal D^\alpha f)(0)
=
\frac12
\int_{\mathbb R}\int_{-\pi}^{\pi}
\mathscr K_\alpha(\tau,\vartheta)
\Bigg[
\int_{\mathbb R}\int_0^{2\pi}
\left|
u(t+\tau,\theta+\vartheta)-u(t,\theta)
\right|
\,d\theta\,dt
\Bigg]
\,d\vartheta\,d\tau.
\end{equation}

For \((\tau,\vartheta)\in\mathbb R\times[-\pi,\pi]\), the fundamental
theorem of calculus yields
\[
u(t+\tau,\theta+\vartheta)-u(t,\theta)
=
\int_0^1
\nabla_{\!c}u(t+r\tau,\theta+r\vartheta)
\cdot(\tau,\vartheta)\,dr.
\]
Therefore,
\[
\int_{\mathbb R}\int_0^{2\pi}
\left|
u(t+\tau,\theta+\vartheta)-u(t,\theta)
\right|
\,d\theta\,dt
\leq
\int_0^1
\int_{\mathbb R}\int_0^{2\pi}
\left|
\nabla_{\!c}u(t+r\tau,\theta+r\vartheta)
\cdot(\tau,\vartheta)
\right|
\,d\theta\,dt\,dr.
\]
For each \(r\in[0,1]\), translation invariance on
\(\mathbb R\times\mathbb S^1\) gives
\[
\int_{\mathbb R}\int_0^{2\pi}
\left|
\nabla_{\!c}u(t+r\tau,\theta+r\vartheta)
\cdot(\tau,\vartheta)
\right|
\,d\theta\,dt
=
\int_{\mathbb R}\int_0^{2\pi}
\left|
\nabla_{\!c}u(t,\theta)
\cdot(\tau,\vartheta)
\right|
\,d\theta\,dt.
\]
Hence
\[
\int_{\mathbb R}\int_0^{2\pi}
\left|
u(t+\tau,\theta+\vartheta)-u(t,\theta)
\right|
\,d\theta\,dt
\leq
\int_{\mathbb R}\int_0^{2\pi}
\left|
\nabla_{\!c}u(t,\theta)\cdot(\tau,\vartheta)
\right|
\,d\theta\,dt.
\]
Together with \eqref{eq:numerator-cylinder-2d}, this gives
\begin{equation}\label{eq:numerator-gradient-bound}
\begin{aligned}
\frac{1}{\gamma_{2,\alpha}}
I_\alpha(\mathcal D^\alpha f)(0)
\leq
\int_{\mathbb R}\int_0^{2\pi}
\Bigg[
\frac12
\int_{\mathbb R}\int_{-\pi}^{\pi}
\left|
\nabla_{\!c}u(t,\theta)\cdot(\tau,\vartheta)
\right|
\mathscr K_\alpha(\tau,\vartheta)
\,d\vartheta\,d\tau
\Bigg]
\,d\theta\,dt.
\end{aligned}
\end{equation}
Setting
\[
E
=
\bigl\{
(t,\theta)\in\mathbb R\times\mathbb S^1:
|\nabla_{\!c}u(t,\theta)|>0
\bigr\},
\]
and observing that the integrand in
\eqref{eq:numerator-gradient-bound} vanishes on \(E^c\), we may rewrite
\eqref{eq:numerator-gradient-bound} as
\begin{equation}\label{eq:numerator-gradient-bound-E}
\begin{aligned}
\frac{1}{\gamma_{2,\alpha}}
I_\alpha(\mathcal D^\alpha f)(0)
\leq
\int_E
\Bigg[
\frac12
\int_{\mathbb R}\int_{-\pi}^{\pi}
\left|
\nabla_{\!c}u(t,\theta)\cdot(\tau,\vartheta)
\right|
\mathscr K_\alpha(\tau,\vartheta)
\,d\vartheta\,d\tau
\Bigg]
\,d\theta\,dt.
\end{aligned}
\end{equation}

For \(\xi=(\xi_1,\xi_2)\in\mathbb S^1\), define
\begin{equation}\label{eq:def-Phi-alpha}
\Phi_\alpha(\xi)
=
\frac12
\int_{\mathbb R}\int_{-\pi}^{\pi}
\left|
\tau\xi_1+\vartheta\xi_2
\right|
\mathscr K_\alpha(\tau,\vartheta)
\,d\vartheta\,d\tau.
\end{equation}
By \eqref{eq:def-K-alpha-2d}, the kernel
\(\mathscr K_\alpha\) is even in each variable. Hence
\[
\Phi_\alpha(\xi_1,\xi_2)
=
\Phi_\alpha(|\xi_1|,|\xi_2|).
\]
For every \(\xi\in\mathbb S^1\), there therefore exists a unique
\(\beta\in[0,\pi/2]\) such that
\[
(|\xi_1|,|\xi_2|)
=
(\cos\beta,\sin\beta).
\]
Using again the evenness of \(\mathscr K_\alpha\)
in each variable, we obtain
\[
\begin{aligned}
\Phi_\alpha(\cos\beta,\sin\beta)
&=
\int_0^\infty\int_0^\pi
\Big(
|\tau\cos\beta+\vartheta\sin\beta|
+
|\tau\cos\beta-\vartheta\sin\beta|
\Big)
\mathscr K_\alpha(\tau,\vartheta)
\,d\vartheta\,d\tau\\
&=
2\int_0^\infty\int_0^\pi
\max\{\tau\cos\beta,\vartheta\sin\beta\}
\mathscr K_\alpha(\tau,\vartheta)
\,d\vartheta\,d\tau\\
&=
A_\alpha(\beta),
\end{aligned}
\]
where we used the elementary identity
\[
|a+b|+|a-b|=2\max\{|a|,|b|\},
\qquad a,b\in\mathbb R.
\]
Consequently,
\begin{equation}\label{eq:sup-Phi-equals-sup-A}
\sup_{\xi\in\mathbb S^1}\Phi_\alpha(\xi)
=
\sup_{0\leq\beta\leq\pi/2}A_\alpha(\beta).
\end{equation}

For every \((t,\theta)\in E\), set
\[
\xi(t,\theta)
=
\frac{
\nabla_{\!c}u(t,\theta)
}{
|\nabla_{\!c}u(t,\theta)|
}
\in\mathbb S^1.
\]
Then, by \eqref{eq:def-Phi-alpha},
\[
\frac12
\int_{\mathbb R}\int_{-\pi}^{\pi}
\left|
\nabla_{\!c}u(t,\theta)\cdot(\tau,\vartheta)
\right|
\mathscr K_\alpha(\tau,\vartheta)
\,d\vartheta\,d\tau
=
|\nabla_{\!c}u(t,\theta)|
\Phi_\alpha(\xi(t,\theta)).
\]
Hence, by \eqref{eq:numerator-gradient-bound-E} and
\eqref{eq:sup-Phi-equals-sup-A},
\[
\begin{aligned}
\frac{1}{\gamma_{2,\alpha}}
I_\alpha(\mathcal D^\alpha f)(0)
&\leq
\left(
\sup_{0\leq\beta\leq\pi/2}A_\alpha(\beta)
\right)
\int_E
|\nabla_{\!c}u(t,\theta)|
\,d\theta\,dt
\\
&=
\left(
\sup_{0\leq\beta\leq\pi/2}A_\alpha(\beta)
\right)
\int_{\mathbb R}\int_0^{2\pi}
|\nabla_{\!c}u(t,\theta)|
\,d\theta\,dt.
\end{aligned}
\]
Together with \eqref{eq:denominator-cylinder-2d},
this gives
\[
R_\alpha(f,0)
\leq
(1-\alpha)
\frac{\gamma_{2,\alpha}}{\gamma_{2,1}}
\sup_{0\leq\beta\leq\pi/2}A_\alpha(\beta).
\]
Taking the supremum over
\(f\in C_c^\infty(\mathbb R^2)\setminus\{0\}\)
yields \eqref{eq:C2-log-polar-upper}.
\end{proof}

\begin{lemma}\label{L2.2}
As \(\alpha\to1^-\),
\begin{equation}\label{eq:uniform-A-alpha-asymptotics}
A_\alpha(\beta)
=
\frac{4}{1-\alpha}+O(1),
\end{equation}
uniformly for \(\beta\in[0,\pi/2]\).
\end{lemma}

\begin{proof}

Throughout the proof, \(C>0\) denotes a constant
independent of \(\alpha\in[1/2,1)\),
\(\beta\in[0,\pi/2]\), and the integration variables,
whose value may change from line to line.

For \(\tau>0\) and \(0<\vartheta<\pi\), set
\[
\rho=(\tau^2+\vartheta^2)^{1/2}.
\]
By \eqref{eq:def-A-alpha-beta}, for every \(\alpha\in[1/2,1)\) and \(\beta\in[0,\pi/2]\), we may write
\begin{equation}\label{eq:A-alpha-singular-decomposition}
A_\alpha(\beta)
=
4\int_{\substack{\tau,\vartheta>0\\
\tau^2+\vartheta^2<1}}
\frac{
\max\{\tau\cos\beta,\vartheta\sin\beta\}
}{
(\tau^2+\vartheta^2)^{1+\alpha/2}
}
\,d\vartheta\,d\tau
+\mathcal R_\alpha(\beta),
\end{equation}
where
\begin{equation}\label{R-alpha-eq}
\begin{aligned}
\mathcal R_\alpha(\beta)
&=
2\int_{\substack{\tau,\vartheta>0\\ \rho<1}}
\max\{\tau\cos\beta,\vartheta\sin\beta\}
\left(
\mathscr K_\alpha(\tau,\vartheta)
-\frac{2}{\rho^{2+\alpha}}
\right)
\,d\vartheta\,d\tau
\\
&\quad+
2\int_{\substack{\tau>0,\ 0<\vartheta<\pi\\ \rho\geq1}}
\max\{\tau\cos\beta,\vartheta\sin\beta\}
\mathscr K_\alpha(\tau,\vartheta)
\,d\vartheta\,d\tau.
\end{aligned}
\end{equation}

We claim that
\begin{equation}\label{eq:A-alpha-remainder-uniform}
\sup_{\substack{
1/2\leq\alpha<1\\
0\leq\beta\leq\pi/2
}}
|\mathcal R_\alpha(\beta)|<\infty.
\end{equation}
Indeed, fix \(\alpha\in[1/2,1)\) and \(\beta\in[0,\pi/2]\). By \eqref{eq:def-K-alpha-2d} and elementary
Taylor estimates, we obtain
\begin{equation}\label{eq:K-alpha-local-expansion}
\left|
\mathscr K_\alpha(\tau,\vartheta)
-
\frac{2}{\rho^{2+\alpha}}
\right|
\leq C\rho^{-\alpha},
\qquad
0<\rho<1.
\end{equation}
Moreover, for every \(\tau>0\) and \(0<\vartheta<\pi\), we have
\[
\max\{\tau\cos\beta,\vartheta\sin\beta\}
\leq \rho.
\]
Thus, by \eqref{eq:K-alpha-local-expansion} and using polar coordinates in the \((\tau,\vartheta)\)-plane, we obtain
\begin{equation}\label{est-R-alpha-rho-less1}
\begin{aligned}
\left|
\int_{\substack{\tau,\vartheta>0\\ \rho<1}}
\max\{\tau\cos\beta,\vartheta\sin\beta\}
\left(
\mathscr K_\alpha(\tau,\vartheta)
-\frac{2}{\rho^{2+\alpha}}
\right)
\,d\vartheta\,d\tau
\right|
& \leq
C\int_0^1\rho^{2-\alpha}\,d\rho\\
& \leq C.
\end{aligned}
\end{equation}
On the complementary region \(\rho\geq1\), the kernel \(\mathscr K_\alpha\) is uniformly bounded for \(0\leq\tau\leq2\), since its denominator is bounded away from zero. Moreover, for \(\tau\geq2\),
\[
\mathscr K_\alpha(\tau,\vartheta)
\leq
Ce^{-\tau/2},
\qquad 0<\vartheta<\pi.
\]
Since
\[
\max\{\tau\cos\beta,\vartheta\sin\beta\}
\leq\tau+\pi,
\qquad \tau>0,\quad 0<\vartheta<\pi,
\]
and
\[
\int_2^\infty(\tau+\pi)e^{-\tau/2}\,d\tau<\infty,
\]
we obtain
\begin{equation}\label{est-R-alpha-rho-less2}
\left|
\int_{\substack{\tau>0,\ 0<\vartheta<\pi\\ \rho\geq1}}
\max\{\tau\cos\beta,\vartheta\sin\beta\}
\mathscr K_\alpha(\tau,\vartheta)
\,d\vartheta\,d\tau
\right|
\leq C.
\end{equation}
Combining \eqref{R-alpha-eq}, \eqref{est-R-alpha-rho-less1},  and \eqref{est-R-alpha-rho-less2}, we obtain \eqref{eq:A-alpha-remainder-uniform}.

We now evaluate the principal term in \eqref{eq:A-alpha-singular-decomposition}. Using the polar change of variables
\[
\tau=\rho\cos\varphi,
\qquad
\vartheta=\rho\sin\varphi,
\qquad
0<\rho<1,\quad 0<\varphi<\frac{\pi}{2},
\]
we obtain for \(\alpha\in [1/2,1)\), 
\[
4\int_{\substack{\tau,\vartheta>0\\
\tau^2+\vartheta^2<1}}
\frac{
\max\{\tau\cos\beta,\vartheta\sin\beta\}
}{
(\tau^2+\vartheta^2)^{1+\alpha/2}
}
\,d\vartheta\,d\tau
=
4\int_0^1\rho^{-\alpha}\,d\rho
\int_0^{\pi/2}
\max\{\cos\varphi\cos\beta,\sin\varphi\sin\beta\}
\,d\varphi.
\]
Splitting the integral at \(\varphi=\pi/2-\beta\) yields
\[
\begin{aligned}
\int_0^{\pi/2}
\max\{\cos\varphi\cos\beta,\sin\varphi\sin\beta\}
\,d\varphi
&=
\cos\beta\sin\left(\frac{\pi}{2}-\beta\right)
+
\sin\beta\cos\left(\frac{\pi}{2}-\beta\right)
\\
&=1.
\end{aligned}
\]
Consequently,
\[
4\int_{\substack{\tau,\vartheta>0\\
\tau^2+\vartheta^2<1}}
\frac{
\max\{\tau\cos\beta,\vartheta\sin\beta\}
}{
(\tau^2+\vartheta^2)^{1+\alpha/2}
}
\,d\vartheta\,d\tau
=
4\int_0^1\rho^{-\alpha}\,d\rho
=
\frac{4}{1-\alpha}.
\]
Combining this identity with
\eqref{eq:A-alpha-singular-decomposition} and
\eqref{eq:A-alpha-remainder-uniform}, we obtain
\eqref{eq:uniform-A-alpha-asymptotics}.
\end{proof}

\begin{lemma}\label{L2.3}
For every \(0<\alpha<1\),
\begin{equation}\label{eq:C2-rad-A-alpha}
C_{2,\alpha}^{\mathrm{rad}}
=
(1-\alpha)\frac{\gamma_{2,\alpha}}{\gamma_{2,1}}
A_\alpha(0).
\end{equation}
Moreover, the supremum defining \(C_{2,\alpha}^{\mathrm{rad}}\) is attained by every nontrivial function \(f\in C_c^\infty(\mathbb R^2)\) that is radially symmetric about the evaluation point and has a monotone radial profile.
\end{lemma}

\begin{proof}
Fix \(0<\alpha<1\). By \eqref{eq:def-A-alpha-beta},
\begin{equation}\label{A-alpha-of-0}
A_\alpha(0)
=
2\int_0^\infty\int_0^\pi
\tau\,\mathscr K_\alpha(\tau,\vartheta)
\,d\vartheta\,d\tau.
\end{equation}

By translation invariance,
it suffices to consider radial functions centered
at the origin. Let
\[
f(x)=F(|x|)\in C_c^\infty(\mathbb R^2)\setminus\{0\},
\]
and set
\[
U(t)=F(e^t),\qquad t\in\mathbb R.
\]

By \eqref{eq:denominator-cylinder-2d},
\begin{equation}\label{eq:radial-denominator-2d}
\frac{I_1(|\nabla f|)(0)}{\gamma_{2,1}}
=
2\pi\int_{\mathbb R}|U'(t)|\,dt.
\end{equation}

Using the evenness of the function
\[
\tau\longmapsto
\int_{\mathbb R}|U(t+\tau)-U(t)|\,dt
\]
and the evenness of \(\mathscr K_\alpha\) in each variable, together with \eqref{eq:numerator-cylinder-2d}, we obtain
\begin{equation}\label{eq:radial-numerator-2d}
\frac{I_\alpha(\mathcal D^\alpha f)(0)}{\gamma_{2,\alpha}}
=
4\pi\int_0^\infty\int_0^\pi
\mathscr K_\alpha(\tau,\vartheta)
\left(
\int_{\mathbb R}
|U(t+\tau)-U(t)|\,dt
\right)
\,d\vartheta\,d\tau.
\end{equation}
Moreover, by Tonelli's theorem, for every \(\tau>0\),
\begin{equation}\label{est-general-U}
\begin{aligned}
\int_{\mathbb R}|U(t+\tau)-U(t)|\,dt
&\leq
\int_{\mathbb R}\int_t^{t+\tau}|U'(s)|\,ds\,dt\\
&=
\tau\int_{\mathbb R}|U'(s)|\,ds.
\end{aligned}
\end{equation}
Combining this estimate with \eqref{A-alpha-of-0},
\eqref{eq:radial-denominator-2d}, and
\eqref{eq:radial-numerator-2d}, we obtain
\[
\begin{aligned}
R_\alpha(f,0)
&\leq
(1-\alpha)\frac{\gamma_{2,\alpha}}{\gamma_{2,1}}
\,2\int_0^\infty\int_0^\pi
\tau\mathscr K_\alpha(\tau,\vartheta)
\,d\vartheta\,d\tau\\
&=
(1-\alpha)\frac{\gamma_{2,\alpha}}{\gamma_{2,1}}
A_\alpha(0).
\end{aligned}
\]
Taking the supremum over the radial class yields
\begin{equation}\label{eq:C2-rad-upper-bound}
C_{2,\alpha}^{\mathrm{rad}}
\leq
(1-\alpha)\frac{\gamma_{2,\alpha}}{\gamma_{2,1}}
A_\alpha(0).
\end{equation}

Assume now that \(F\) is monotone. Then \(U\) is monotone and, for every \(\tau>0\) and \(t\in\mathbb R\),
\[
|U(t+\tau)-U(t)|
=
\int_t^{t+\tau}|U'(s)|\,ds.
\]
Consequently, equality holds in \eqref{est-general-U}, and hence
\[
R_\alpha(f,0)
=
(1-\alpha)\frac{\gamma_{2,\alpha}}{\gamma_{2,1}}
A_\alpha(0).
\]
Thus, the upper bound in \eqref{eq:C2-rad-upper-bound} is attained by every such function. This proves \eqref{eq:C2-rad-A-alpha} and shows that the supremum is a maximum.
\end{proof}

We next compute \(A_\alpha(0)\). Recall that, for
\(a\in\mathbb R\) and \(k\in\mathbb N_0\), the
Pochhammer symbol is defined by
\[
(a)_k
=
\begin{cases}
1, & k=0,\\[1mm]
a(a+1)\cdots(a+k-1), & k\geq1.
\end{cases}
\]
For \(a\notin\{0,-1,-2,\ldots\}\), we also have
\[
(a)_k
=
\frac{\Gamma(a+k)}{\Gamma(a)}.
\]

\begin{lemma}\label{L2.4}
For every \(0<\alpha<1\), we have
\begin{equation}\label{eq:explicit-A-alpha-0}
A_\alpha(0)
=
\frac{2\pi}{\alpha^2}
\left(
\frac{2\Gamma(1-\alpha)}
{\Gamma(1-\alpha/2)^2}
-1
\right).
\end{equation}
\end{lemma}

\begin{proof}
Fix \(0<\alpha<1\). By \eqref{A-alpha-of-0},
\[
A_\alpha(0)
=
2\int_0^\infty\int_0^\pi
\tau\,\mathscr K_\alpha(\tau,\vartheta)
\,d\vartheta\,d\tau.
\]
Using \eqref{eq:def-K-alpha-2d} and
\[
2(\cosh\tau-\cos\vartheta)
=
e^\tau
\bigl(1-2e^{-\tau}\cos\vartheta+e^{-2\tau}\bigr),
\]
we obtain
\[
\mathscr K_\alpha(\tau,\vartheta)
=
\frac{
e^{-\alpha\tau}+e^{-2\tau}
}{
\bigl(1-2e^{-\tau}\cos\vartheta+e^{-2\tau}\bigr)^{1+\alpha/2}
}.
\]
Hence
\begin{equation}\label{eq:A-alpha-0-integral}
A_\alpha(0)
=
2\int_0^\infty
\tau\bigl(e^{-\alpha\tau}+e^{-2\tau}\bigr)
\left[
\int_0^\pi
\frac{d\vartheta}
{\bigl(1-2e^{-\tau}\cos\vartheta+e^{-2\tau}\bigr)^{1+\alpha/2}}
\right]
\,d\tau.
\end{equation}

Set
\[
a=\frac{\alpha}{2}\in\left(0,\frac12\right),
\qquad
r=e^{-\tau}\in(0,1).
\]
Then
\[
\int_0^\pi
\frac{d\vartheta}
{\bigl(1-2e^{-\tau}\cos\vartheta+e^{-2\tau}\bigr)^{1+\alpha/2}}
=
\int_0^\pi
\frac{d\vartheta}
{\bigl(1-2r\cos\vartheta+r^2\bigr)^{1+a}}.
\]
Using the standard hypergeometric identity
(see, e.g., \cite[Section~2.1]{AndrewsAskeyRoy})
\[
\int_0^\pi
\frac{d\vartheta}
{\bigl(1-2r\cos\vartheta+r^2\bigr)^{1+a}}
=
\pi
\sum_{k=0}^\infty
\left(
\frac{(1+a)_k}{k!}
\right)^2
r^{2k},
\qquad 0<r<1,
\]
we obtain
\[
\int_0^\pi
\frac{d\vartheta}
{\bigl(1-2e^{-\tau}\cos\vartheta+e^{-2\tau}\bigr)^{1+\alpha/2}}
=
\pi
\sum_{k=0}^\infty
\left(
\frac{(1+a)_k}{k!}
\right)^2
e^{-2k\tau}.
\]
Combining this identity with
\eqref{eq:A-alpha-0-integral}, applying Tonelli's theorem, and using
\[
\int_0^\infty \tau e^{-\lambda\tau}\,d\tau
=
\frac{1}{\lambda^2},
\qquad \lambda>0,
\]
we obtain
\[
A_\alpha(0)
=
\frac{\pi}{2}
\sum_{k=0}^\infty
\left(
\frac{(1+a)_k}{k!}
\right)^2
\left[
\frac{1}{(k+a)^2}
+
\frac{1}{(k+1)^2}
\right].
\]
Using
\[
\frac{(1+a)_k}{k+a}
=
\frac{(a)_k}{a},
\qquad
\frac{(1+a)_k}{k+1}
=
\frac{(a)_{k+1}}{a(k+1)},
\]
we obtain
\[
A_\alpha(0)
=
\frac{\pi}{2a^2}
\left[
2\sum_{k=0}^\infty
\left(\frac{(a)_k}{k!}\right)^2
-1
\right].
\]
Since \(0<a<1/2\), Gauss' summation formula
(see, e.g., \cite[Theorem~2.2.2]{AndrewsAskeyRoy}) gives
\[
\sum_{k=0}^\infty
\left(\frac{(a)_k}{k!}\right)^2
=
{}_2F_1(a,a;1;1)
=
\frac{\Gamma(1-2a)}
{\Gamma(1-a)^2}.
\]
Substituting this identity and recalling that \(a=\alpha/2\), 
we obtain \eqref{eq:explicit-A-alpha-0}.

\end{proof}

We now prove Theorems~\ref{thm:leading-asymptotics-2d} and~\ref{thm:sharp-constant-radial-2d}.

\subsection*{Proof of Theorem~\ref{thm:leading-asymptotics-2d}}

\begin{proof}
By Lemmas~\ref{L2.1} and~\ref{L2.3}, we obtain, for every \(0<\alpha<1\),
\[
(1-\alpha)\frac{\gamma_{2,\alpha}}{\gamma_{2,1}}
A_\alpha(0)
=
C_{2,\alpha}^{\mathrm{rad}}
\leq C_{2,\alpha}^{*}
\leq
(1-\alpha)\frac{\gamma_{2,\alpha}}{\gamma_{2,1}}
\sup_{\beta\in[0,\pi/2]}A_\alpha(\beta).
\]
By Lemma~\ref{L2.2},
\[
A_\alpha(\beta)
=
\frac{4}{1-\alpha}+O(1),
\qquad \alpha\to1^-,
\]
uniformly for \(\beta\in[0,\pi/2]\).

On the other hand, the normalization \eqref{eq:gamma-n-beta} yields
\[
\frac{\gamma_{2,\alpha}}{\gamma_{2,1}}
=
2^{1-\alpha}
\frac{\Gamma(1-\alpha/2)}{\Gamma(\alpha/2)},
\qquad 0<\alpha<1.
\]
By the smoothness of the Gamma function,
\[
\frac{\gamma_{2,\alpha}}{\gamma_{2,1}}
=
1+O(1-\alpha),
\qquad \alpha\to1^-.
\]
Combining these estimates with the preceding two-sided inequality, we obtain \eqref{eq:leading-asymptotics-C2}.
\end{proof}

\subsection*{Proof of Theorem~\ref{thm:sharp-constant-radial-2d}}

\begin{proof}
By Lemmas~\ref{L2.3} and~\ref{L2.4}, and recalling \eqref{eq:gamma-n-beta}, we obtain, for every \(0<\alpha<1\),
\[
\begin{aligned}
C_{2,\alpha}^{\mathrm{rad}}
&=
(1-\alpha)\frac{\gamma_{2,\alpha}}{\gamma_{2,1}}A_\alpha(0)\\
&=
\frac{2^{2-\alpha}\pi(1-\alpha)}{\alpha^2}
\frac{\Gamma(1-\alpha/2)}{\Gamma(\alpha/2)}
\left(
2\frac{\Gamma(1-\alpha)}
{\Gamma(1-\alpha/2)^2}-1
\right).
\end{aligned}
\]
This establishes \eqref{eq:sharp-C2-rad}. The assertion concerning the attainment of the supremum follows directly from Lemma~\ref{L2.3}.
\end{proof}

\section{The sharp constant \(C_{n,\alpha}^{*}\) for \(n\geq3\)}
\label{sec3}

This section is devoted to the proof of
Theorem~\ref{thm:sharp-constant-higher-dim}
and the derivation of Corollary~\ref{cor:strict-HMP}.
Throughout this section, we assume that \(n\geq3\).
We begin with several auxiliary lemmas.

We first derive logarithmic-polar representations
of the Riesz potentials involved in our problem.
For \(0<\alpha<1\), \(\tau\in\mathbb R\), and
\(0\leq\vartheta\leq\pi\), with
\((\tau,\vartheta)\neq(0,0)\), define
\begin{equation}\label{eq:def-K-alpha-n}
\mathscr K_{n,\alpha}(\tau,\vartheta)
=
\frac{
2\cosh\!\left(\frac{n-\alpha}{2}\tau\right)
}{
\bigl[2(\cosh\tau-\cos\vartheta)\bigr]^{(n+\alpha)/2}
}.
\end{equation}

\begin{lemma}
Let \(0<\alpha<1\) and \(f\in C_c^\infty(\mathbb R^n)\). Set
\[
u(t,\omega)=f(e^t\omega),
\qquad
(t,\omega)\in\mathbb R\times\mathbb S^{n-1}.
\]
Then
\begin{equation}\label{eq:denominator-cylinder-n}
\frac{I_1(|\nabla f|)(0)}{\gamma_{n,1}}
=
\int_{\mathbb R}\int_{\mathbb S^{n-1}}
|\nabla_c u(t,\omega)|
\,d\sigma(\omega)\,dt,
\end{equation}
where \(\nabla_c\) denotes the gradient on the cylinder \(\mathbb R\times\mathbb S^{n-1}\), defined by
\[
\nabla_c u(t,\omega)
=
\bigl(\partial_t u(t,\omega),
\nabla_{\mathbb S^{n-1}}u(t,\omega)\bigr),
\qquad
(t,\omega)\in\mathbb R\times\mathbb S^{n-1},
\]
so that
\[
|\nabla_c u(t,\omega)|^2
=
|\partial_t u(t,\omega)|^2
+
|\nabla_{\mathbb S^{n-1}}u(t,\omega)|^2.
\]
Moreover,
\begin{equation}\label{eq:numerator-cylinder-n}
\begin{aligned}
&\frac{I_\alpha(\mathcal D^\alpha f)(0)}
{\gamma_{n,\alpha}}\\
&\quad=
\frac12
\int_{\mathbb R}\int_{\mathbb R}
\int_{\mathbb S^{n-1}}\int_{\mathbb S^{n-1}}
\mathscr K_{n,\alpha}
\bigl(\tau,\arccos(\omega\cdot\nu)\bigr)
|u(t+\tau,\nu)-u(t,\omega)|
\,d\sigma(\nu)\,d\sigma(\omega)\,d\tau\,dt.
\end{aligned}
\end{equation}
\end{lemma}

\begin{proof}
We first establish \eqref{eq:denominator-cylinder-n}.  For every \(t\in\mathbb R\) and \(\omega\in\mathbb S^{n-1}\), the definition of \(u\) gives
\[
\partial_t u(t,\omega)
=
e^t\nabla f(e^t\omega)\cdot\omega
\]
and
\[
\nabla_{\mathbb S^{n-1}}u(t,\omega)
=
e^t\left(
\nabla f(e^t\omega)
-
\bigl(\nabla f(e^t\omega)\cdot\omega\bigr)\omega
\right).
\]	
Consequently,
\[
\begin{aligned}
|\nabla_c u(t,\omega)|^2
&=
|\partial_t u(t,\omega)|^2
+
|\nabla_{\mathbb S^{n-1}}u(t,\omega)|^2\\
&=
e^{2t}|\nabla f(e^t\omega)|^2,
\end{aligned}
\]
and hence
\begin{equation}\label{eq:cylinder-gradient-n}
|\nabla_c u(t,\omega)|
=
e^t|\nabla f(e^t\omega)|,
\qquad
(t,\omega)\in\mathbb R\times\mathbb S^{n-1}.
\end{equation}
Using the polar change of variables
\[
y=e^t\omega,
\qquad
(t,\omega)\in\mathbb R\times\mathbb S^{n-1},
\]
we obtain
\[
\begin{aligned}
\frac{I_1(|\nabla f|)(0)}{\gamma_{n,1}}
&=
\int_{\mathbb R^n}
\frac{|\nabla f(y)|}{|y|^{n-1}}\,dy\\
&=
\int_{\mathbb R}\int_{\mathbb S^{n-1}}
e^t|\nabla f(e^t\omega)|
\,d\sigma(\omega)\,dt\\
&=
\int_{\mathbb R}\int_{\mathbb S^{n-1}}
|\nabla_c u(t,\omega)|
\,d\sigma(\omega)\,dt,
\end{aligned}
\]
where the last equality follows from \eqref{eq:cylinder-gradient-n}. This proves \eqref{eq:denominator-cylinder-n}.

We now prove \eqref{eq:numerator-cylinder-n}.  By the definitions of \(I_\alpha\) and \(\mathcal D^\alpha\),
\[
\frac{I_\alpha(\mathcal D^\alpha f)(0)}
{\gamma_{n,\alpha}}
=
\int_{\mathbb R^n}\int_{\mathbb R^n}
\frac{|f(y)-f(z)|}
{|y|^{n-\alpha}|y-z|^{n+\alpha}}
\,dz\,dy.
\]
Using the changes of variables
\[
y=e^t\omega,\qquad z=e^s\nu,
\]
where \(t,s\in\mathbb R\) and \(\omega,\nu\in\mathbb S^{n-1}\), we have
\[
dy=e^{nt}\,dt\,d\sigma(\omega),
\qquad
dz=e^{ns}\,ds\,d\sigma(\nu),
\]
and
\[
|y-z|^2
=
2e^{t+s}\bigl(\cosh(s-t)-\omega\cdot\nu\bigr).
\]
Consequently,
\[
\frac{I_\alpha(\mathcal D^\alpha f)(0)}
{\gamma_{n,\alpha}}
=
\int_{\mathbb R}\int_{\mathbb R}
\int_{\mathbb S^{n-1}}\int_{\mathbb S^{n-1}}
\frac{
e^{(n-\alpha)(s-t)/2}
|u(s,\nu)-u(t,\omega)|
}{
[2(\cosh(s-t)-\omega\cdot\nu)]^{(n+\alpha)/2}
}
\,d\sigma(\nu)\,d\sigma(\omega)\,ds\,dt.
\]
Making the change of variables
\[
(t,s)\mapsto(t,\tau),
\qquad \tau=s-t,
\]
we obtain
\begin{equation}\label{eq:numerator-cylinder-n-intermediate}
\frac{I_\alpha(\mathcal D^\alpha f)(0)}
{\gamma_{n,\alpha}}
=
\int_{\mathbb R}\int_{\mathbb R}
\int_{\mathbb S^{n-1}}\int_{\mathbb S^{n-1}}
\frac{
e^{(n-\alpha)\tau/2}
|u(t+\tau,\nu)-u(t,\omega)|
}{
[2(\cosh\tau-\omega\cdot\nu)]^{(n+\alpha)/2}
}
\,d\sigma(\nu)\,d\sigma(\omega)\,d\tau\,dt.
\end{equation}
Making the change of variables
\[
(t,\tau,\omega,\nu)
\mapsto
(t+\tau,-\tau,\nu,\omega)
\]
in \eqref{eq:numerator-cylinder-n-intermediate}, we also obtain
\[
\frac{I_\alpha(\mathcal D^\alpha f)(0)}
{\gamma_{n,\alpha}}
=
\int_{\mathbb R}\int_{\mathbb R}
\int_{\mathbb S^{n-1}}\int_{\mathbb S^{n-1}}
\frac{
e^{-(n-\alpha)\tau/2}
|u(t+\tau,\nu)-u(t,\omega)|
}{
[2(\cosh\tau-\omega\cdot\nu)]^{(n+\alpha)/2}
}
\,d\sigma(\nu)\,d\sigma(\omega)\,d\tau\,dt.
\]
Averaging these two representations yields
\[
\begin{aligned}
&\frac{I_\alpha(\mathcal D^\alpha f)(0)}
{\gamma_{n,\alpha}}
\\
& \quad =
\int_{\mathbb R}\int_{\mathbb R}
\int_{\mathbb S^{n-1}}\int_{\mathbb S^{n-1}}
\frac{
\cosh\!\left(\frac{n-\alpha}{2}\tau\right)
}{
[2(\cosh\tau-\omega\cdot\nu)]^{(n+\alpha)/2}
}|u(t+\tau,\nu)-u(t,\omega)|
\,d\sigma(\nu)\,d\sigma(\omega)\,d\tau\,dt.
\end{aligned}
\]
Finally, by \eqref{eq:def-K-alpha-n},
we obtain \eqref{eq:numerator-cylinder-n}.
This completes the proof.
\end{proof}

\begin{lemma}\label{L3.2}
Let \(0<\alpha<1\) and
\[
f(x)=F(|x|)\in C_c^\infty(\mathbb R^n)\setminus\{0\}.
\]
Then
\begin{equation}\label{eq:radial-bound-n}
R_\alpha(f,0)
\leq
(1-\alpha)
\frac{\gamma_{n,\alpha}}{\gamma_{n,1}}
\mathfrak C_{n,\alpha}^{\rm rad}.
\end{equation}
Moreover, equality holds whenever \(F\) is monotone.
\end{lemma}

\begin{proof}
We first establish \eqref{eq:radial-bound-n}.
Set
\[
U(t)=F(e^t),\qquad t\in\mathbb R.
\]
Since \(u(t,\omega)=U(t)\) for every
\((t,\omega)\in\mathbb R\times\mathbb S^{n-1}\),
we have
\[
\partial_t u(t,\omega)=U'(t),
\qquad
\nabla_{\mathbb S^{n-1}}u(t,\omega)=0.
\]
Consequently, \eqref{eq:denominator-cylinder-n} gives
\begin{equation}\label{eq:radial-denominator-n}
\frac{I_1(|\nabla f|)(0)}{\gamma_{n,1}}
=
|\mathbb S^{n-1}|
\int_{\mathbb R}|U'(t)|\,dt.
\end{equation}

On the other hand, for every \(\omega\in\mathbb S^{n-1}\), using spherical coordinates with \(\theta=\arccos(\omega\cdot\nu)\), we obtain
\[
\int_{\mathbb S^{n-1}}
\mathscr K_{n,\alpha}
\bigl(\tau,\arccos(\omega\cdot\nu)\bigr)
\,d\sigma(\nu)
=
|\mathbb S^{n-2}|
\int_0^\pi
\mathscr K_{n,\alpha}(\tau,\theta)
(\sin\theta)^{n-2}\,d\theta,
\qquad \tau\in\mathbb R\setminus\{0\}.
\]
Using this identity, the evenness of
\(\mathscr K_{n,\alpha}\) in \(\tau\), and
the evenness of the function
\[
\tau\longmapsto
\int_{\mathbb R}|U(t+\tau)-U(t)|\,dt,
\]
we deduce from \eqref{eq:numerator-cylinder-n} that
\begin{equation}\label{eq:radial-numerator-n}
\begin{aligned}
&\frac{I_\alpha(\mathcal D^\alpha f)(0)}
{\gamma_{n,\alpha}}
\\
&\quad =
|\mathbb S^{n-1}||\mathbb S^{n-2}|
\int_0^\infty\int_0^\pi
\mathscr K_{n,\alpha}(\tau,\theta)
(\sin\theta)^{n-2}
\left(
\int_{\mathbb R}
|U(t+\tau)-U(t)|\,dt
\right)
\,d\theta\,d\tau.
\end{aligned}
\end{equation}
Moreover, by the fundamental theorem of calculus
and Tonelli's theorem, for every \(\tau>0\),
\[
\int_{\mathbb R}|U(t+\tau)-U(t)|\,dt
\leq
\tau\int_{\mathbb R}|U'(s)|\,ds.
\]
Consequently, \eqref{eq:def-radial-coefficient-n} and \eqref{eq:radial-numerator-n}
 yield
\begin{equation}\label{eq-gen-inequality}
\frac{I_\alpha(\mathcal D^\alpha f)(0)}
{\gamma_{n,\alpha}}
\leq
|\mathbb S^{n-1}|
\mathfrak C_{n,\alpha}^{\rm rad}
\int_{\mathbb R}|U'(t)|\,dt.
\end{equation}
Combining \eqref{eq:radial-denominator-n} and \eqref{eq-gen-inequality}, we obtain \eqref{eq:radial-bound-n}.

Assume now that \(F\) is monotone. Then \(U\) is monotone and, for every \(\tau>0\) and \(t\in\mathbb R\),
\[
|U(t+\tau)-U(t)|
=
\int_t^{t+\tau}|U'(s)|\,ds.
\]
Consequently,
\[
\int_{\mathbb R}|U(t+\tau)-U(t)|\,dt
=
\tau\int_{\mathbb R}|U'(s)|\,ds.
\]
Thus, equality holds in \eqref{eq-gen-inequality}, and hence
\[
R_\alpha(f,0)
=
(1-\alpha)
\frac{\gamma_{n,\alpha}}{\gamma_{n,1}}
\mathfrak C_{n,\alpha}^{\rm rad}.
\]
This completes the proof.
\end{proof}

For \(0<\alpha<1\) and \(\beta\in[0,\pi/2]\), define
\begin{equation}\label{eq:def-A-n-alpha-beta}
A_{n,\alpha}(\beta)
=
\int_0^\infty\int_0^\pi
\mathscr K_{n,\alpha}(\tau,\theta)
(\sin\theta)^{n-2}\int_{\mathbb S^{n-2}}
\max\left\{
\tau\cos\beta,\,
\theta\sin\beta\,|\eta\cdot e|
\right\}
\,d\sigma(\eta)\,d\theta\,d\tau,
\end{equation}
where \(e\in\mathbb S^{n-2}\) is any fixed unit vector. By rotational invariance, \(A_{n,\alpha}(\beta)\) does not depend on the choice of \(e\).

\begin{lemma}\label{L3.3}
For every \(0<\alpha<1\) and every
\(f\in C_c^\infty(\mathbb R^n)\setminus\{0\}\),
\begin{equation}\label{eq:Cn-angular-upper}
R_\alpha(f,0)
\leq
(1-\alpha)
\frac{\gamma_{n,\alpha}}{\gamma_{n,1}}
\sup_{\beta\in[0,\pi/2]}A_{n,\alpha}(\beta).
\end{equation}
Moreover,
\begin{equation}\label{eq:An-radial-value}
A_{n,\alpha}(0)
=
\mathfrak C_{n,\alpha}^{\rm rad}.
\end{equation}
\end{lemma}

\begin{proof}

Fix \(0<\alpha<1\) and \(f\in C_c^\infty(\mathbb R^n)\setminus\{0\}\).
Set
\[
u(t,\omega)=f(e^t\omega),
\qquad
(t,\omega)\in\mathbb R\times\mathbb S^{n-1}.
\]
For each \(\omega\in\mathbb S^{n-1}\), let
\[
\mathbb S_\omega^{n-2}
=
\{\eta\in\mathbb S^{n-1}:\eta\cdot\omega=0\},
\]
and denote by \(d\sigma_\omega\) the induced
surface measure on \(\mathbb S_\omega^{n-2}\).
Using spherical coordinates centered at \(\omega\),
we write
\[
\nu=\cos\theta\,\omega+\sin\theta\,\eta,
\qquad
0<\theta<\pi,\quad
\eta\in\mathbb S_\omega^{n-2}.
\]
Then
\[
d\sigma(\nu)
=
(\sin\theta)^{n-2}
\,d\sigma_\omega(\eta)\,d\theta.
\]
By \eqref{eq:numerator-cylinder-n}, the evenness
of \(\mathscr K_{n,\alpha}\) in its first variable,
and the symmetry obtained by interchanging the
two spherical variables and translating the
\(t\)-variable, we have
\begin{equation}\label{eq:L33-numerator}
\frac{I_\alpha(\mathcal D^\alpha f)(0)}
{\gamma_{n,\alpha}}
=
\int_0^\infty\int_0^\pi
\mathscr K_{n,\alpha}(\tau,\theta)
(\sin\theta)^{n-2}
\mathcal J_u(\tau,\theta)
\,d\theta\,d\tau,
\end{equation}
where
\[
\mathcal J_u(\tau,\theta)
=
\int_{\mathbb R}\int_{\mathbb S^{n-1}}
\int_{\mathbb S_\omega^{n-2}}
\left|
u(t+\tau,\cos\theta\,\omega+\sin\theta\,\eta)
-u(t,\omega)
\right|
\,d\sigma_\omega(\eta)\,d\sigma(\omega)\,dt.
\]

Fix \(\tau>0\), \(0<\theta<\pi\),
\(\omega\in\mathbb S^{n-1}\), and
\(\eta\in\mathbb S_\omega^{n-2}\).
For \(s\in[0,1]\), define
\[
\omega_s
=
\cos(s\theta)\omega+\sin(s\theta)\eta
\]
and
\[
\eta_s
=
-\sin(s\theta)\omega+\cos(s\theta)\eta.
\]
Then
\[
\eta_s\in\mathbb S_{\omega_s}^{n-2},
\qquad
\frac{d\omega_s}{ds}=\theta\eta_s.
\]
By the fundamental theorem of calculus,
\[
u(t+\tau,\omega_1)-u(t,\omega)
=
\int_0^1
\Big[
\tau\,\partial_tu(t+s\tau,\omega_s)
+
\theta\,
\nabla_{\mathbb S^{n-1}}u(t+s\tau,\omega_s)
\cdot\eta_s
\Big]\,ds.
\]
Consequently,
\begin{equation}\label{diff-u-eq}
|u(t+\tau,\omega_1)-u(t,\omega)|
\leq
\int_0^1
\left|
\tau\,\partial_tu(t+s\tau,\omega_s)
+
\theta\,
\nabla_{\mathbb S^{n-1}}u(t+s\tau,\omega_s)
\cdot\eta_s
\right|\,ds.
\end{equation}

For each \(s\in[0,1]\), the transformation
\[
(\omega,\eta)\longmapsto(\omega_s,\eta_s)
\]
is the geodesic flow at time \(s\theta\) on the
unit tangent bundle of \(\mathbb S^{n-1}\).
In particular, it preserves the natural measure
\[
d\sigma(\omega)\,d\sigma_\omega(\eta).
\]
Combining this measure invariance with translation
invariance of Lebesgue measure in the \(t\)-variable
and using \eqref{diff-u-eq}, we obtain
\begin{equation}\label{eq:L33-J-bound}
\mathcal J_u(\tau,\theta)
\leq
\int_{\mathbb R}\int_{\mathbb S^{n-1}}
\int_{\mathbb S_\omega^{n-2}}
\left|
\tau\,\partial_tu(t,\omega)
+
\theta\,\nabla_{\mathbb S^{n-1}}u(t,\omega)
\cdot\eta
\right|
\,d\sigma_\omega(\eta)\,d\sigma(\omega)\,dt.
\end{equation}

Next, using the symmetry \(\eta\mapsto-\eta\)
and the elementary identity
\[
\frac{|a+b|+|a-b|}{2}
=
\max\{|a|,|b|\},
\qquad a,b\in\mathbb R,
\]
we deduce that
\begin{equation}\label{eq:L33-symmetry}
\begin{aligned}
&\int_{\mathbb S_\omega^{n-2}}
\left|
\tau\,\partial_tu(t,\omega)
+
\theta\,\nabla_{\mathbb S^{n-1}}u(t,\omega)
\cdot\eta
\right|
\,d\sigma_\omega(\eta)\\
&\quad =
\int_{\mathbb S_\omega^{n-2}}
\max\left\{
\tau|\partial_tu(t,\omega)|,\,
\theta
\left|
\nabla_{\mathbb S^{n-1}}u(t,\omega)
\cdot\eta
\right|
\right\}
\,d\sigma_\omega(\eta).
\end{aligned}
\end{equation}

Set
\[
G(t,\omega)=|\nabla_cu(t,\omega)|.
\]
At every point where \(G(t,\omega)>0\),
choose \(\beta(t,\omega)\in[0,\pi/2]\) such that
\[
\cos\beta(t,\omega)
=
\frac{|\partial_tu(t,\omega)|}{G(t,\omega)}
\]
and
\[
\sin\beta(t,\omega)
=
\frac{
|\nabla_{\mathbb S^{n-1}}u(t,\omega)|
}{G(t,\omega)}.
\]
At points where \(G(t,\omega)=0\), set
\(\beta(t,\omega)=0\).
By rotational invariance of the surface measure
on \(\mathbb S_\omega^{n-2}\),
\eqref{eq:L33-symmetry} becomes
\begin{equation}\label{eq:L33-angular}
\begin{aligned}
&\int_{\mathbb S_\omega^{n-2}}
\left|
\tau\,\partial_tu(t,\omega)
+
\theta\,\nabla_{\mathbb S^{n-1}}u(t,\omega)
\cdot\eta
\right|
\,d\sigma_\omega(\eta)
\\
&\quad=
G(t,\omega)
\int_{\mathbb S^{n-2}}
\max\left\{
\tau\cos\beta(t,\omega),\,
\theta\sin\beta(t,\omega)|\eta\cdot e|
\right\}
\,d\sigma(\eta),
\end{aligned}
\end{equation}
where \(e\in\mathbb S^{n-2}\) is any fixed
unit vector.

Combining \eqref{eq:L33-numerator},
\eqref{eq:L33-J-bound}, and
\eqref{eq:L33-angular}, and applying
Tonelli's theorem, we obtain
\[
\begin{aligned}
\frac{I_\alpha(\mathcal D^\alpha f)(0)}
{\gamma_{n,\alpha}}
&\leq
\int_{\mathbb R}\int_{\mathbb S^{n-1}}
G(t,\omega)
A_{n,\alpha}(\beta(t,\omega))
\,d\sigma(\omega)\,dt
\\
&\leq
\left(
\sup_{\beta\in[0,\pi/2]}A_{n,\alpha}(\beta)
\right)
\int_{\mathbb R}\int_{\mathbb S^{n-1}}
|\nabla_cu(t,\omega)|
\,d\sigma(\omega)\,dt.
\end{aligned}
\]
By \eqref{eq:denominator-cylinder-n},
\[
\int_{\mathbb R}\int_{\mathbb S^{n-1}}
|\nabla_cu(t,\omega)|
\,d\sigma(\omega)\,dt
=
\frac{I_1(|\nabla f|)(0)}{\gamma_{n,1}}.
\]
Hence,
\[
I_\alpha(\mathcal D^\alpha f)(0)
\leq
\frac{\gamma_{n,\alpha}}{\gamma_{n,1}}
\left(
\sup_{\beta\in[0,\pi/2]}A_{n,\alpha}(\beta)
\right)
I_1(|\nabla f|)(0).
\]
Since \(f\not\equiv0\), we have
\(I_1(|\nabla f|)(0)>0\).
Multiplying by \(1-\alpha\) and dividing by
\(I_1(|\nabla f|)(0)\), we obtain
\eqref{eq:Cn-angular-upper}.

Finally, taking \(\beta=0\) in
\eqref{eq:def-A-n-alpha-beta}, we find
\[
\begin{aligned}
A_{n,\alpha}(0)
&=
\int_0^\infty\int_0^\pi
\mathscr K_{n,\alpha}(\tau,\theta)
(\sin\theta)^{n-2}
\int_{\mathbb S^{n-2}}\tau
\,d\sigma(\eta)\,d\theta\,d\tau
\\
&=
|\mathbb S^{n-2}|
\int_0^\infty\int_0^\pi
\tau\,\mathscr K_{n,\alpha}(\tau,\theta)
(\sin\theta)^{n-2}
\,d\theta\,d\tau
\\
&=
\mathfrak C_{n,\alpha}^{\rm rad},
\end{aligned}
\]
where the last equality follows from
\eqref{eq:def-radial-coefficient-n}.
This proves \eqref{eq:An-radial-value}
and completes the proof.
\end{proof}

\begin{lemma}\label{L3.4}
For every \(0<\alpha<1\) and every
\(\beta\in[0,\pi/2]\),
\begin{equation}\label{eq:An-max-radial}
A_{n,\alpha}(\beta)
\leq A_{n,\alpha}(0).
\end{equation}
Consequently,
\[
\sup_{\beta\in[0,\pi/2]}A_{n,\alpha}(\beta)
=
\mathfrak C_{n,\alpha}^{\rm rad}.
\]
\end{lemma}

\begin{proof}
Fix \(0<\alpha<1\) and \(0\leq \beta\leq \pi/2\). Define
\[
\rho(\tau,\theta)
=
\mathscr K_{n,\alpha}(\tau,\theta)
\left(\frac{\sin\theta}{\theta}\right)^{n-2},
\qquad
\tau\in\mathbb R,\quad 0<\theta<\pi.
\]
We extend \(\rho\) by zero for \(\theta\geq\pi\).

By rotational invariance, we may assume that
\[
e=(1,0,\ldots,0)\in\mathbb R^{n-1}.
\]
For \(0<\theta<\pi\) and
\(\eta\in\mathbb S^{n-2}\), set
\[
y=\theta\eta\in\mathbb R^{n-1}.
\]
Then
\[
|y|=\theta,\qquad
y_1=\theta(\eta\cdot e),
\qquad
dy=\theta^{n-2}\,d\theta\,d\sigma(\eta).
\]
Consequently, by \eqref{eq:def-A-n-alpha-beta}
and the definition of \(\rho\), we obtain
\[
A_{n,\alpha}(\beta)
=
\int_0^\infty\int_{|y|<\pi}
\max\left\{
\tau\cos\beta,\,
|y_1|\sin\beta
\right\}
\rho(\tau,|y|)
\,dy\,d\tau.
\]
Since \(\rho(\tau,\theta)\) has been extended
by zero for \(\theta\geq\pi\), we may equivalently
write
\[
A_{n,\alpha}(\beta)
=
\int_0^\infty\int_{\mathbb R^{n-1}}
\max\left\{
\tau\cos\beta,\,
|y_1|\sin\beta
\right\}
\rho(\tau,|y|)
\,dy\,d\tau.
\]
Using the elementary identity
\[
\max\{a,|b|\}
=
\frac{|a+b|+|a-b|}{2},
\qquad a\geq0,\quad b\in\mathbb R,
\]
and the evenness of \(\rho\) in \(\tau\),
we obtain
\begin{equation}\label{simpl-A-n-alpha-beta}
\begin{aligned}
A_{n,\alpha}(\beta)
&=
\frac12\int_0^\infty\int_{\mathbb R^{n-1}}
\Big(
|\tau\cos\beta+y_1\sin\beta|
+
|\tau\cos\beta-y_1\sin\beta|
\Big)
\rho(\tau,|y|)
\,dy\,d\tau\\
&=
\frac12
\int_{\mathbb R}\int_{\mathbb R^{n-1}}
|\tau\cos\beta+y_1\sin\beta|
\rho(\tau,|y|)
\,dy\,d\tau.
\end{aligned}
\end{equation}

Set
\[
z=(\tau,y)\in\mathbb R\times \mathbb R^{n-1},
\qquad
v_\beta=(\cos\beta,\sin\beta,0,\ldots,0)
\in\mathbb S^{n-1}.
\]
For every \(r>0\), define
\[
F_r(s)
=
\rho\left(rs,r\sqrt{1-s^2}\right),
\qquad -1<s<1,
\]
with endpoint values understood by continuity.
Since \(\rho\) is even in its first variable,
the function \(F_r\) is even.

Using polar coordinates \(z=r\xi\), where
\(r>0\) and \(\xi\in\mathbb S^{n-1}\), we have
\[
\tau=r\xi_1,
\qquad
|y|=r\sqrt{1-\xi_1^2},
\qquad
dz=r^{n-1}\,dr\,d\sigma(\xi).
\]
Consequently, by \eqref{simpl-A-n-alpha-beta},
we obtain
\begin{equation}\label{eq:L34-spherical-representation}
A_{n,\alpha}(\beta)
=
\frac12\int_0^\infty r^n
\int_{\mathbb S^{n-1}}
|\xi\cdot v_\beta|
F_r(|\xi_1|)
\,d\sigma(\xi)\,dr.
\end{equation}

We now show that, for every \(r>0\),
the function \(F_r\) is nondecreasing on
\([0,1]\). Fix \(r>0\). For \(0<s<1\), set
\[
\theta=r\sqrt{1-s^2}.
\]
Then
\[
rs=\sqrt{r^2-\theta^2}.
\]
Since
\[
\frac{d\theta}{ds}
=
-\frac{rs}{\sqrt{1-s^2}}<0
\]
and \(\rho\) has a continuous zero extension
for \(\theta\geq\pi\), it suffices to prove
that the function
\[
\theta\longmapsto
\rho\left(\sqrt{r^2-\theta^2},\theta\right)
\]
is decreasing on
\(0<\theta<\min\{r,\pi\}\).

 Set
\[
a=\frac{n-\alpha}{2},
\qquad
p=a+\alpha.
\]
For
\(0<\theta<\min\{r,\pi\}\), let
\[
\tau=\sqrt{r^2-\theta^2}.
\]
Recalling the definitions of \(\rho\) and
\(\mathscr K_{n,\alpha}\), and observing that
\[
\frac{d\tau}{d\theta}=-\frac{\theta}{\tau},
\]
direct differentiation gives
\[
\frac{d}{d\theta}
\log\rho\left(\sqrt{r^2-\theta^2},\theta\right)
=
-\frac{a\theta}{\tau}\tanh(a\tau)
+
p\frac{\theta\sinh\tau/\tau-\sin\theta}
{\cosh\tau-\cos\theta}
+
(n-2)\left(\cot\theta-\frac1\theta\right).
\]
Using the  elementary inequalities
\[
\frac{\theta\sinh\tau/\tau-\sin\theta}
{\cosh\tau-\cos\theta}
\leq
\frac{\theta}{\tau}\tanh\frac{\tau}{2},
\qquad
\frac1\theta-\cot\theta
\geq\frac{\theta}{3},
\]
we obtain 
\begin{equation}\label{eq:L34-rho-upper}
\frac{d}{d\theta}
\log\rho\left(\sqrt{r^2-\theta^2},\theta\right)
\leq
\frac{\theta}{\tau}D_{n,\alpha}(\tau),
\end{equation}
where
\[
D_{n,\alpha}(\tau)
=
p\tanh\frac{\tau}{2}
-a\tanh(a\tau)
-\frac{n-2}{3}\tau.
\]

Using the elementary inequalities
\[
2\tanh\frac{\tau}{2}-\tanh \tau<\frac{\tau}{3},\qquad \tanh\frac{\tau}{2}-\tanh \tau<0,\qquad 
\tau>0,
\]
the monotonicity of \(\tanh\), \(a>1\), and \(p=a+\alpha\), we obtain for every \(\tau>0\), 
\[
\begin{aligned}
D_{n,\alpha}(\tau)
&=a\left[\tanh\frac{\tau}{2}-\tanh(a\tau)\right]
+\alpha\tanh\frac{\tau}{2}
-\frac{n-2}{3}\tau\\
&\leq
a\left[\tanh\frac{\tau}{2}-\tanh\tau\right]
+\tanh\frac{\tau}{2}-\frac{\tau}{3}\\
&\leq
2\tanh\frac{\tau}{2}-\tanh\tau-\frac{\tau}{3}\\
&<0. 
\end{aligned}
\]
Since \(\rho(\tau,\theta)>0\) for
\(\tau>0\) and \(0<\theta<\pi\),
it follows from \eqref{eq:L34-rho-upper} that
\[
\frac{d}{d\theta}
\rho\left(\sqrt{r^2-\theta^2},\theta\right)
<0,
\qquad
0<\theta<\min\{r,\pi\}.
\]
Consequently, \(F_r\) is nondecreasing on \([0,1]\).

We next prove that, for every \(r>0\),
\begin{equation}\label{eq:L34-spherical-rearrangement}
\int_{\mathbb S^{n-1}}
|\xi\cdot v_\beta|F_r(|\xi_1|)
\,d\sigma(\xi)
\leq
\int_{\mathbb S^{n-1}}
|\xi_1|F_r(|\xi_1|)
\,d\sigma(\xi).
\end{equation}
To this end, for almost every
\(\xi\in\mathbb S^{n-1}\), write
\[
\xi=(q\cos\varphi,q\sin\varphi,\zeta),
\]
where
\[
\zeta\in B^{n-2},
\qquad
q=\sqrt{1-|\zeta|^2},
\qquad
0\leq\varphi<2\pi,
\]
and \(B^{n-2}\) denotes the unit ball in
\(\mathbb R^{n-2}\). In these coordinates,
the surface measure satisfies
\[
d\sigma(\xi)=d\varphi\,d\zeta.
\]
Moreover,
\[
|\xi\cdot v_\beta|
=
q|\cos(\varphi-\beta)|,
\qquad
|\xi_1|=q|\cos\varphi|.
\]
Consequently,
\begin{equation}\label{int-S-n-1-eq}
\int_{\mathbb S^{n-1}}
|\xi\cdot v_\beta|F_r(|\xi_1|)
\,d\sigma(\xi)
=
\int_{B^{n-2}}q
\int_0^{2\pi}
|\cos(\varphi-\beta)|
F_r(q|\cos\varphi|)
\,d\varphi\,d\zeta.
\end{equation}

For each fixed \(q\in[0,1]\), the functions
\[
\varphi\longmapsto F_r(q|\cos\varphi|)
\]
and
\[
\varphi\longmapsto|\cos\varphi|
\]
are nonnegative, even, and \(\pi\)-periodic
on \(\mathbb R\). Moreover, since \(F_r\)
is nondecreasing on \([0,1]\), both functions
are nonincreasing on \([0,\pi/2]\). Hence, by the circular rearrangement inequality
(see \cite{BaernsteinCircle}), we obtain
\[
\int_0^{2\pi}
F_r(q|\cos\varphi|)
|\cos(\varphi-\beta)|\,d\varphi
\leq
\int_0^{2\pi}
F_r(q|\cos\varphi|)
|\cos\varphi|\,d\varphi.
\]
Together with \eqref{int-S-n-1-eq}, this gives
\[
\begin{aligned}
\int_{\mathbb S^{n-1}}
|\xi\cdot v_\beta|F_r(|\xi_1|)
\,d\sigma(\xi) &\leq \int_{B^{n-2}}	\int_0^{2\pi}q|\cos\varphi|F_r(q|\cos\varphi|)
\,d\varphi\,d\zeta\\
&=\int_{\mathbb S^{n-1}}
|\xi_1|F_r(|\xi_1|)
\,d\sigma(\xi),
\end{aligned}
\]
which proves \eqref{eq:L34-spherical-rearrangement}. 

Finally, combining
\eqref{eq:L34-spherical-representation}
and \eqref{eq:L34-spherical-rearrangement}, we obtain
\[
\begin{aligned}
A_{n,\alpha}(\beta)
&\leq
\frac12\int_0^\infty r^n
\int_{\mathbb S^{n-1}}
|\xi_1|F_r(|\xi_1|)
\,d\sigma(\xi)\,dr
\\
&=
A_{n,\alpha}(0).
\end{aligned}
\]
This proves \eqref{eq:An-max-radial}. The final assertion follows from
\eqref{eq:An-radial-value}.
\end{proof}

We are now ready to prove
Theorem~\ref{thm:sharp-constant-higher-dim}.

\subsection*{Proof of Theorem~\ref{thm:sharp-constant-higher-dim}}

\begin{proof}
By Lemmas~\ref{L3.3} and~\ref{L3.4}, for every
\(0<\alpha<1\) and every
\(f\in C_c^\infty(\mathbb R^n)\setminus\{0\}\),
we have
\[
\begin{aligned}
R_\alpha(f,0)
&\leq
(1-\alpha)
\frac{\gamma_{n,\alpha}}{\gamma_{n,1}}
\sup_{\beta\in[0,\pi/2]}A_{n,\alpha}(\beta)\\
&=
(1-\alpha)
\frac{\gamma_{n,\alpha}}{\gamma_{n,1}}
\mathfrak C_{n,\alpha}^{\rm rad}.
\end{aligned}
\]
By translation invariance, taking the supremum
over all admissible functions and evaluation
points yields
\[
C_{n,\alpha}^{*}
\leq
(1-\alpha)
\frac{\gamma_{n,\alpha}}{\gamma_{n,1}}
\mathfrak C_{n,\alpha}^{\rm rad}.
\]
On the other hand, Lemma~\ref{L3.2} shows
that equality holds for every nontrivial
smooth radially monotone function centered
at the evaluation point. Consequently,
\[
C_{n,\alpha}^{*}
=
(1-\alpha)
\frac{\gamma_{n,\alpha}}{\gamma_{n,1}}
\mathfrak C_{n,\alpha}^{\rm rad}.
\]
This proves \eqref{eq:sharp-Cn} and the
assertion concerning the attainment of
the optimal constant.
\end{proof}

\subsection*{Proof of Corollary~\ref{cor:strict-HMP}}

\begin{proof}

Fix \(0<\alpha<1\). Let
\[
f(y)=F(|y|)\in C_c^\infty(\mathbb R^n)
\setminus\{0\},
\]
where \(F\) is monotone. By Theorem~\ref{thm:sharp-constant-higher-dim},
\begin{equation}\label{eq:cor-HMP-extremizer}
R_\alpha(f,0)=C_{n,\alpha}^{*}.
\end{equation}

Choose \(R>0\) such that
\(\operatorname{supp}f\subset B(0,R)\),
where \(B(0,R)\) denotes the open ball
of radius \(R\) centered at the origin
in \(\mathbb R^n\). Since \(F\) is monotone and nontrivial,
we have \(f(0)\neq0\). For every \(y\in\mathbb R^n\) with
\(|y|>R\), taking \(h=-2y\) gives
\[
f(y+h)=f(y)=0,
\]
whereas
\[
|h|\int_0^1|\nabla f(y+th)|\,dt>0,
\]
since the segment joining \(y\) to \(-y\)
passes through the origin, where
\(f(0)\neq0\). Hence, by continuity, there exists
a neighborhood \(V_{-2y}\) of \(-2y\)
such that
\[
|f(y+h)-f(y)|
<
|h|\int_0^1|\nabla f(y+th)|\,dt,
\qquad h\in V_{-2y}.
\]
On the other hand, the fundamental theorem
of calculus gives
\[
|f(y+h)-f(y)|
\leq
|h|\int_0^1|\nabla f(y+th)|\,dt
\]
for every \(h\in\mathbb R^n\). Consequently, for  \(|y|>R\), splitting the integral
defining \(\mathcal D^\alpha f(y)\)
over \(V_{-2y}\) and its complement,
we obtain
\[
\begin{aligned}
\mathcal D^\alpha f(y)
&=
\int_{V_{-2y}}
\frac{|f(y+h)-f(y)|}{|h|^{n+\alpha}}\,dh
+
\int_{\mathbb R^n\setminus V_{-2y}}
\frac{|f(y+h)-f(y)|}{|h|^{n+\alpha}}\,dh\\
&<
\int_{\mathbb R^n}
\frac{1}{|h|^{n+\alpha-1}}
\int_0^1|\nabla f(y+th)|\,dt\,dh\\
&=
\int_0^1t^{\alpha-1}\,dt
\int_{\mathbb R^n}
\frac{|\nabla f(z)|}{|y-z|^{n+\alpha-1}}\,dz\\
&=
\frac{1}{\alpha\gamma_{n,1-\alpha}}
I_{1-\alpha}(|\nabla f|)(y),
\end{aligned}
\]
that is,
\begin{equation}\label{str-ineq-eq}
\mathcal D^\alpha f(y)
<
\frac{1}{\alpha\gamma_{n,1-\alpha}}
I_{1-\alpha}(|\nabla f|)(y),
\qquad |y|>R.
\end{equation}
Moreover, the same argument yields
\begin{equation}\label{large-ineq-eq}
\mathcal D^\alpha f(y)
\leq
\frac{1}{\alpha\gamma_{n,1-\alpha}}
I_{1-\alpha}(|\nabla f|)(y),
\qquad y\in\mathbb R^n.
\end{equation}

Since \eqref{str-ineq-eq} holds on the
set \(\{|y|>R\}\) of positive measure,
while \eqref{large-ineq-eq} holds on
\(\mathbb R^n\), the strict positivity
of the Riesz kernel yields
\[
I_\alpha(\mathcal D^\alpha f)(0)
<
\frac{1}{\alpha\gamma_{n,1-\alpha}}
I_\alpha I_{1-\alpha}(|\nabla f|)(0).
\]
By the semigroup property
\(I_\alpha I_{1-\alpha}=I_1\), we obtain
\[
\begin{aligned}
(1-\alpha)I_\alpha(\mathcal D^\alpha f)(0) &<
\frac{1-\alpha}{\alpha\gamma_{n,1-\alpha}}
I_1(|\nabla f|)(0)\\
&=
C_{n,\alpha}^{\mathrm{HMP}}
I_1(|\nabla f|)(0),
\end{aligned}
\]
where the last equality follows from \eqref{eq:gamma-n-beta} and \eqref{eq-C-HMP}.  Since \(I_1(|\nabla f|)(0)>0\),
\eqref{eq:cor-HMP-extremizer} yields
\[
C_{n,\alpha}^{*}
<
C_{n,\alpha}^{\mathrm{HMP}}.
\]

Finally, recalling that
\[
\lim_{\alpha\to1^-}C_{n,\alpha}^{*}=K_n,
\qquad
\lim_{\alpha\to1^-}
C_{n,\alpha}^{\mathrm{HMP}}
=|\mathbb S^{n-1}|,
\]
we conclude that
\[
\lim_{\alpha\to1^-}
\frac{C_{n,\alpha}^{*}}
{C_{n,\alpha}^{\mathrm{HMP}}}
=
\frac{K_n}{|\mathbb S^{n-1}|}<1.
\]
This completes the proof.
\end{proof}

\section{First-order asymptotics of the normalized Riesz potential}
\label{sec4}

This section is devoted to the proof of
Theorem~\ref{thm:first-order-expansion}.
Throughout this section, we assume that \(n\geq2\).
We begin with two auxiliary lemmas.

\begin{lemma}\label{L4.1}
Let \(f\in C_c^2(\mathbb R^n)\) and \(p\in\{1,\infty\}\). Then, as
\(\alpha\to1^-\),
\begin{equation}\label{eq:D-alpha-first-order}
\left\|
(1-\alpha)\mathcal D^\alpha f
-
K_n|\nabla f|
-
(1-\alpha)\mathcal F_f
\right\|_{L^p(\mathbb R^n)}
=
O\bigl((1-\alpha)^2\bigr),
\end{equation}
where \(\mathcal F_f\) is defined by
\eqref{eq:def-Ff}.
\end{lemma}

\begin{proof}
Fix  \(0<\alpha<1\). For \(x,h\in\mathbb R^n\), set
\[
Q_f(x,h)
=
|f(x+h)-f(x)|
-
|\nabla f(x)\cdot h|.
\]
By the definition of \(\mathcal D^\alpha f\),
for every \(x\in \mathbb{R}^n\), 
\[
\mathcal D^\alpha f(x)
=
\int_{|h|<1}
\frac{|\nabla f(x)\cdot h|}
{|h|^{n+\alpha}}\,dh
+
\int_{|h|<1}
\frac{Q_f(x,h)}
{|h|^{n+\alpha}}\,dh
+
\int_{|h|\geq1}
\frac{|f(x+h)-f(x)|}
{|h|^{n+\alpha}}\,dh.
\]
Using polar coordinates and the definition
of \(K_n\), we obtain
\[
\begin{aligned}
\int_{|h|<1}
\frac{|\nabla f(x)\cdot h|}
{|h|^{n+\alpha}}\,dh
&=
K_n|\nabla f(x)|
\int_0^1r^{-\alpha}\,dr
\\
&=
\frac{K_n}{1-\alpha}|\nabla f(x)|.
\end{aligned}
\]
Consequently, by \eqref{eq:def-Ff},
\begin{equation}\label{eq:L41-remainder-identity}
\begin{aligned}
&(1-\alpha)\mathcal D^\alpha f(x)
-
K_n|\nabla f(x)|
-
(1-\alpha)\mathcal F_f(x)
\\
&\quad=
(1-\alpha)
\int_{|h|<1}Q_f(x,h)
\left(
\frac1{|h|^{n+\alpha}}
-
\frac1{|h|^{n+1}}
\right)\,dh
\\
&\qquad+
(1-\alpha)
\int_{|h|\geq1}|f(x+h)-f(x)|
\left(
\frac1{|h|^{n+\alpha}}
-
\frac1{|h|^{n+1}}
\right)\,dh.
\end{aligned}
\end{equation}
We now estimate the two integrals on the
right-hand side. 

By Taylor's theorem, for every \(x,h\in \mathbb{R}^n\), 
\[
|Q_f(x,h)|
\leq
|h|^2\int_0^1(1-s)
|D^2f(x+sh)|\,ds,
\]
where \(|D^2f|\) denotes the operator norm of the Hessian matrix of \(f\).
Hence, taking the supremum for \(p=\infty\),
and applying Tonelli's theorem and translation
invariance for \(p=1\), we obtain, for every
\(h\in\mathbb R^n\),
\begin{equation}\label{eq:L41-Q-bound}
\|Q_f(\cdot,h)\|_{L^p(\mathbb R^n)}
\leq
\frac{|h|^2}{2}
\|D^2f\|_{L^p(\mathbb R^n)}.
\end{equation}
Moreover,
\begin{equation}\label{eq:L41-difference-bound}
\|f(\cdot+h)-f(\cdot)\|_{L^p(\mathbb R^n)}
\leq
2\|f\|_{L^p(\mathbb R^n)}.
\end{equation}

On the other hand, for \(0<r<1\),
\[
\begin{aligned}
\left|
\frac1{r^{n+\alpha}}
-
\frac1{r^{n+1}}
\right|
&=
\frac{1-r^{1-\alpha}}{r^{n+1}}
\\
&\leq
(1-\alpha)\frac{|\log r|}{r^{n+1}}.
\end{aligned}
\]
Similarly, for \(r\geq1\),
\[
\begin{aligned}
\left|
\frac1{r^{n+\alpha}}
-
\frac1{r^{n+1}}
\right|
&=
\frac{1-r^{-(1-\alpha)}}{r^{n+\alpha}}
\\
&\leq
(1-\alpha)\frac{\log r}{r^{n+\alpha}}.
\end{aligned}
\]
Combining these estimates with \eqref{eq:L41-remainder-identity},
\eqref{eq:L41-Q-bound}, and \eqref{eq:L41-difference-bound},
taking the \(L^p\)-norm, and using
polar coordinates, we obtain
\[
\left\|
(1-\alpha)\mathcal D^\alpha f
-
K_n|\nabla f|
-
(1-\alpha)\mathcal F_f
\right\|_{L^p(\mathbb R^n)}
\leq
(1-\alpha)^2|\mathbb S^{n-1}|
\left(
\frac12\|D^2f\|_{L^p(\mathbb R^n)}
+
\frac{2}{\alpha^2}
\|f\|_{L^p(\mathbb R^n)}
\right).
\]
Since \(f\in C_c^2(\mathbb R^n)\) and \(\alpha^{-2}\) remains bounded as
\(\alpha\to1^-\), we obtain \eqref{eq:D-alpha-first-order}. This completes the proof.
\end{proof}

\begin{lemma}\label{L4.2}
Let \(g\in L^1(\mathbb R^n)\cap L^\infty(\mathbb R^n)\). Then, as \(\alpha\to1^-\),
\begin{equation}\label{eq:I-alpha-first-order}
I_\alpha g(x)
=
I_1g(x)
-
(1-\alpha)\gamma_{n,1}
\int_{\mathbb R^n}
\frac{\log|x-y|+\lambda_n}{|x-y|^{n-1}}
g(y)\,dy
+
O\bigl((1-\alpha)^2\bigr),
\end{equation}
uniformly for \(x\in\mathbb R^n\), where
\(\lambda_n\) is defined by
\eqref{eq:def-lambda-n}.
Moreover, there exists a constant \(C_n>0\),
depending only on \(n\), such that
\begin{equation}\label{ineq-normI-alpha-g}
\|I_\alpha g\|_{L^\infty(\mathbb R^n)}
\leq
C_n\left(
\|g\|_{L^1(\mathbb R^n)}
+
\|g\|_{L^\infty(\mathbb R^n)}
\right)
\end{equation}
for every \(\alpha\in[1/2,1]\).
\end{lemma}

\begin{proof}

Throughout the proof, \(C_n>0\) denotes a constant
depending only on \(n\), whose value may change
from line to line.

For \(1/2\leq s\leq1\) and
\(z\in\mathbb R^n\setminus\{0\}\), set
\[
k_s(z)=\gamma_{n,s}|z|^{s-n}.
\]
By \eqref{eq:gamma-n-beta}, we have
\[
\frac{d}{ds}\log\gamma_{n,s}
=
-\frac12\psi\left(\frac{n-s}{2}\right)
-\log2
-\frac12\psi\left(\frac{s}{2}\right).
\]
In particular, by \eqref{eq:def-lambda-n},
\[
\left.
\frac{d}{ds}\log\gamma_{n,s}
\right|_{s=1}
=
\lambda_n.
\]
Consequently,
\begin{equation}\label{eq:L42-kernel-derivative}
\left.\frac{\partial k_s(z)}{\partial s}
\right|_{s=1}
=
\gamma_{n,1}
\frac{\log|z|+\lambda_n}{|z|^{n-1}}.
\end{equation}

We now estimate the second derivative of
\(k_s\) with respect to \(s\). Set
\[
b_n(s)=\frac{d}{ds}\log\gamma_{n,s}.
\]
Direct differentiation gives
\[
\frac{\partial^2k_s(z)}{\partial s^2}
=
k_s(z)
\left[
\bigl(\log|z|+b_n(s)\bigr)^2
+b_n'(s)
\right].
\]
Since \(s\mapsto\gamma_{n,s}\) is positive and smooth
on \([1/2,1]\), the functions
\(\gamma_{n,s}\), \(b_n(s)\), and \(b_n'(s)\)
are bounded on this interval. Consequently, there exists
\(C_n>0\) such that
\begin{equation}\label{est-sderv-ks}
\left|
\frac{\partial^2k_s(z)}{\partial s^2}
\right|
\leq
C_n
\begin{cases}
|z|^{1/2-n}(1+\log^2|z|),
&0<|z|<1,\\
|z|^{1-n}(1+\log^2|z|),
&|z|\geq1,
\end{cases}
\end{equation}
uniformly for \(s\in[1/2,1]\).

For \(1/2\leq\alpha<1\) and
\(z\in\mathbb R^n\setminus\{0\}\), set
\[
\mathcal R_{n,\alpha}(z)
=
k_\alpha(z)-k_1(z)
+(1-\alpha)
\left.\frac{\partial k_s(z)}{\partial s}
\right|_{s=1}.
\]
Taylor's theorem gives
\[
|\mathcal R_{n,\alpha}(z)|
\leq
\frac{(1-\alpha)^2}{2}
\sup_{\alpha\leq s\leq1}
\left|
\frac{\partial^2k_s(z)}{\partial s^2}
\right|.
\]
By \eqref{est-sderv-ks}, it follows that
\begin{equation}\label{est-R-alpha-z}
|\mathcal R_{n,\alpha}(z)|
\leq
C_n(1-\alpha)^2
\begin{cases}
|z|^{1/2-n}(1+\log^2|z|),
&0<|z|<1,\\
|z|^{1-n}(1+\log^2|z|),
&|z|\geq1.
\end{cases}
\end{equation}

Using the definitions of \(I_\alpha\) and \(I_1\),
together with \eqref{eq:L42-kernel-derivative},
we have, for \(1/2\leq\alpha<1\) and
\(x\in\mathbb R^n\),
\[
\begin{aligned}
I_\alpha g(x)-I_1g(x)
&=
\int_{\mathbb R^n}
(k_\alpha(x-y)-k_1(x-y))g(y)\,dy\\
&=
\int_{\mathbb R^n}
(k_\alpha(z)-k_1(z))g(x-z)\,dz\\
&=
\int_{\mathbb R^n}
\mathcal R_{n,\alpha}(z)g(x-z)\,dz-
(1-\alpha)
\int_{\mathbb R^n}
\left.\frac{\partial k_s(z)}{\partial s}
\right|_{s=1}g(x-z)\,dz\\
&=
\int_{\mathbb R^n}
\mathcal R_{n,\alpha}(z)g(x-z)\,dz
-(1-\alpha)\gamma_{n,1}
\int_{\mathbb R^n}
\frac{\log|x-y|+\lambda_n}
{|x-y|^{n-1}}g(y)\,dy.
\end{aligned}
\]
Hence,
\[
I_\alpha g(x)-I_1g(x)
+
(1-\alpha)\gamma_{n,1}
\int_{\mathbb R^n}
\frac{\log|x-y|+\lambda_n}
{|x-y|^{n-1}}g(y)\,dy
=
\int_{\mathbb R^n}
\mathcal R_{n,\alpha}(z)g(x-z)\,dz.
\]
Together with \eqref{est-R-alpha-z}, this yields
\[
\begin{aligned}
&\Bigg|
I_\alpha g(x)-I_1g(x)
+(1-\alpha)\gamma_{n,1}
\int_{\mathbb R^n}
\frac{\log|x-y|+\lambda_n}
{|x-y|^{n-1}}g(y)\,dy
\Bigg|
\\
&\quad\leq
C_n(1-\alpha)^2
\Bigg[
\|g\|_{L^\infty(\mathbb R^n)}
\int_{|z|<1}
\frac{1+\log^2|z|}{|z|^{n-1/2}}
\,dz
+
\|g\|_{L^1(\mathbb R^n)}
\sup_{|z|\geq1}
\frac{1+\log^2|z|}{|z|^{n-1}}
\Bigg].
\end{aligned}
\]
Since \(n\geq2\), we have
\[
\int_{|z|<1}
\frac{1+\log^2|z|}{|z|^{n-1/2}}
\,dz
=
|\mathbb S^{n-1}|
\int_0^1r^{-1/2}(1+\log^2 r)\,dr
<\infty
\]
and
\[
\sup_{|z|\geq1}
\frac{1+\log^2|z|}{|z|^{n-1}}
<\infty.
\]
Therefore,
\[
\Bigg\|
I_\alpha g-I_1g
+
(1-\alpha)\gamma_{n,1}
\int_{\mathbb R^n}
\frac{\log|\cdot-y|+\lambda_n}
{|\cdot-y|^{n-1}}g(y)\,dy
\Bigg\|_{L^\infty(\mathbb R^n)}
\leq
C_n(1-\alpha)^2
\left(
\|g\|_{L^1(\mathbb R^n)}
+
\|g\|_{L^\infty(\mathbb R^n)}
\right).
\]
This proves \eqref{eq:I-alpha-first-order}.

We now prove the second assertion. For
\(1/2\leq\alpha\leq1\) and \(x\in\mathbb R^n\),
we have
\[
\begin{aligned}
|I_\alpha g(x)|
&\leq
\gamma_{n,\alpha}
\int_{\mathbb R^n}
\frac{|g(y)|}{|x-y|^{n-\alpha}}\,dy\\
&\leq
\gamma_{n,\alpha}
\left(
\|g\|_{L^\infty(\mathbb R^n)}
\int_{|z|<1}\frac{dz}{|z|^{n-\alpha}}
+
\|g\|_{L^1(\mathbb R^n)}
\right)\\
&=
\gamma_{n,\alpha}
\left(
\frac{|\mathbb S^{n-1}|}{\alpha}
\|g\|_{L^\infty(\mathbb R^n)}
+
\|g\|_{L^1(\mathbb R^n)}
\right).
\end{aligned}
\]
Since \(\gamma_{n,\alpha}\) is bounded on
\([1/2,1]\), we obtain \eqref{ineq-normI-alpha-g}.
This completes the proof.
\end{proof}

\subsection*{Proof of Theorem~\ref{thm:first-order-expansion}}

\begin{proof}
Let \(f\in C_c^2(\mathbb R^n)\) and set
\[
\varepsilon=1-\alpha.
\]

We first note that
\[
\mathcal F_f\in L^1(\mathbb R^n)\cap
L^\infty(\mathbb R^n).
\] 
Indeed, by \eqref{eq:def-Ff},
\eqref{eq:L41-Q-bound}, and
\eqref{eq:L41-difference-bound}, for
\(p\in\{1,\infty\}\), we have
\[
\begin{aligned}
\|\mathcal F_f\|_{L^p(\mathbb R^n)}
&\leq
\frac12\|D^2f\|_{L^p(\mathbb R^n)}
\int_{|h|<1}\frac{dh}{|h|^{n-1}}
+
2\|f\|_{L^p(\mathbb R^n)}
\int_{|h|\geq1}\frac{dh}{|h|^{n+1}}\\
&=
|\mathbb S^{n-1}|
\left(
\frac12\|D^2f\|_{L^p(\mathbb R^n)}
+
2\|f\|_{L^p(\mathbb R^n)}
\right)\\
&<\infty.
\end{aligned}
\]

For \(g\in L^1(\mathbb R^n)\cap
L^\infty(\mathbb R^n)\), define
\[
\mathcal L_n g(x)
=
\gamma_{n,1}
\int_{\mathbb R^n}
\frac{\log|x-y|+\lambda_n}
{|x-y|^{n-1}}g(y)\,dy.
\]
Since \(n\geq2\),
\[
\int_{|z|<1}
\frac{1+|\log|z||}{|z|^{n-1}}\,dz
<\infty,
\qquad
\sup_{|z|\geq1}
\frac{1+\log|z|}{|z|^{n-1}}
<\infty.
\]
Consequently, there exists a constant \(C_n>0\),
depending only on \(n\), such that
\begin{equation}\label{eq:log-potential-bound}
\|\mathcal L_n g\|_{L^\infty(\mathbb R^n)}
\leq
C_n\left(
\|g\|_{L^1(\mathbb R^n)}
+
\|g\|_{L^\infty(\mathbb R^n)}
\right).
\end{equation}

By Lemma~\ref{L4.1}, we can write
\[
\varepsilon\mathcal D^\alpha f
=
K_n|\nabla f|
+
\varepsilon\mathcal F_f
+
E_\alpha,
\]
where
\begin{equation}\label{eq:proof-thm-Ealpha}
\|E_\alpha\|_{L^1(\mathbb R^n)}
+
\|E_\alpha\|_{L^\infty(\mathbb R^n)}
=
O(\varepsilon^2),\qquad \varepsilon\to 0^+. 
\end{equation}
Applying \(I_\alpha\) and using its
linearity, we obtain
\begin{equation}\label{eq:proof-thm-decomposition}
\begin{aligned}
\varepsilon I_\alpha(\mathcal D^\alpha f)
&=
K_n I_\alpha(|\nabla f|)
+
\varepsilon I_\alpha(\mathcal F_f)
+
I_\alpha E_\alpha.
\end{aligned}
\end{equation}
We now estimate the three terms on the right-hand side.

By the first assertion of Lemma~\ref{L4.2},
\begin{equation}\label{eq:proof-thm-gradient}
I_\alpha(|\nabla f|)
=
I_1(|\nabla f|)
-
\varepsilon\mathcal L_n(|\nabla f|)
+
O_{L^\infty}(\varepsilon^2). 
\end{equation}
Applying the same lemma to
\(g=\mathcal F_f\) and using
\eqref{eq:log-potential-bound}, gives
\begin{equation}\label{eq:proof-thm-Ff}
I_\alpha(\mathcal F_f)
=
I_1(\mathcal F_f)
+
O_{L^\infty}(\varepsilon). 
\end{equation}
Finally, the second assertion of
Lemma~\ref{L4.2}, together with
\eqref{eq:proof-thm-Ealpha}, yields
\begin{equation}\label{eq:proof-thm-remainder}
\|I_\alpha E_\alpha\|_{L^\infty(\mathbb R^n)}
\leq
C_n\left(
\|E_\alpha\|_{L^1(\mathbb R^n)}
+
\|E_\alpha\|_{L^\infty(\mathbb R^n)}
\right)
=
O(\varepsilon^2). 
\end{equation}

Substituting \eqref{eq:proof-thm-gradient},
\eqref{eq:proof-thm-Ff}, and
\eqref{eq:proof-thm-remainder} into
\eqref{eq:proof-thm-decomposition}, we obtain
\[
\varepsilon I_\alpha(\mathcal D^\alpha f)
=
K_n I_1(|\nabla f|)
+
\varepsilon\left(
I_1(\mathcal F_f)
-
K_n\mathcal L_n(|\nabla f|)
\right)
+
O_{L^\infty}(\varepsilon^2).
\]
By \eqref{eq:def-Af},
\[
I_1(\mathcal F_f)
-
K_n\mathcal L_n(|\nabla f|)
=
\mathcal A_f.
\]
Therefore, \eqref{eq:first-order-expansion} follows.
\end{proof}

\section{Almost-everywhere convergence for \(W^{1,1}\) functions}
\label{sec5}

In this section, we prove the almost-everywhere
convergence stated in Theorem~\ref{thm:W11-full-limit}
for functions in \(W^{1,1}(\mathbb R^n)\),
without passing to a subsequence.
Throughout this section, we assume that \(n\ge 2\).

We first establish the following lemma.

\begin{lemma}\label{L5.1}
There exists a constant \(C_n>0\), depending
only on \(n\), with the following property:
for every \(u,v\in W^{1,1}(\mathbb R^n)\),
there exists a set \(E_{u,v}\subset\mathbb R^n\)
of full measure such that
\begin{equation}\label{eq:uniform-stability-W11}
\bigl|
(1-\alpha)I_\alpha(\mathcal D^\alpha u)(x)
-
(1-\alpha)I_\alpha(\mathcal D^\alpha v)(x)
\bigr|
\leq
C_n I_1(|\nabla(u-v)|)(x)
<\infty
\end{equation}
for every \(x\in E_{u,v}\) and every
\(\alpha\in[1/2,1)\).
\end{lemma}

\begin{proof}

By \eqref{eq:HMP-comparison} and the fact that
\(C_{n,\alpha}^{\mathrm{HMP}}\) is uniformly
bounded for \(\alpha\in[1/2,1)\), there exists
a constant \(C_n>0\), depending only on \(n\),
such that
\begin{equation}\label{eq:L51-smooth-HMP}
(1-\alpha)I_\alpha(\mathcal D^\alpha\varphi)(x)
\leq
C_n I_1(|\nabla\varphi|)(x)
\end{equation}
for every \(\varphi\in C_c^\infty(\mathbb R^n)\), every \(x\in\mathbb R^n\), and every
\(\alpha\in[1/2,1)\).

We claim that for every \(w\in W^{1,1}(\mathbb R^n)\), there exists a set \(E_w\subset \mathbb R^n\) of full measure such that 
\begin{equation}\label{claim-L5.1}
(1-\alpha)I_\alpha(\mathcal D^\alpha w)(x)
\leq
C_n I_1(|\nabla w|)(x)	<\infty 
\end{equation}
for every \(x\in E_w\) and every \(\alpha\in[1/2,1)\). 

Indeed, fix \(w\in W^{1,1}(\mathbb R^n)\).  By the density of \(C_c^\infty(\mathbb R^n)\)
in \(W^{1,1}(\mathbb R^n)\), we can choose
\((w_j)\subset C_c^\infty(\mathbb R^n)\)
such that
\begin{equation}\label{cv-W11}
\|w_j-w\|_{W^{1,1}(\mathbb R^n)}
\leq 2^{-j},
\qquad j\geq1.
\end{equation}
Consequently,
\[
G=\sum_{j=1}^\infty
|\nabla(w_j-w)|
\in L^1(\mathbb R^n).
\]
Since the Riesz potential \(I_1h\) is finite
almost everywhere whenever
\(h\in L^1(\mathbb R^n)\), the set
\[
E_w=
\left\{
x\in\mathbb R^n:
I_1(G+|\nabla w|)(x)<\infty
\right\}
\]
has full measure. Moreover, for every \(x\in E_w\), Tonelli's theorem gives
\[
\sum_{j=1}^{\infty}
I_1(|\nabla(w_j-w)|)(x)
=
I_1G(x)<\infty.
\]
Hence, for every \(x\in E_w\),
\begin{equation}\label{lim-gradients}
\left|
I_1(|\nabla w_j|)(x)
-
I_1(|\nabla w|)(x)
\right|
\leq
I_1(|\nabla(w_j-w)|)(x)
\longrightarrow0
\end{equation}
as \(j\to\infty\).

Now fix \(x\in E_w\) and \(\alpha\in[1/2,1)\). By \eqref{cv-W11}, \(w_j\to w\) almost everywhere
in \(\mathbb R^n\) as \(j\to\infty\).  
Hence, applying Fatou's lemma to the
double integral defining \(I_\alpha(\mathcal D^\alpha w)(x)\), and using \eqref{eq:L51-smooth-HMP}
and \eqref{lim-gradients}, we obtain
\[
\begin{aligned}
(1-\alpha)I_\alpha(\mathcal D^\alpha w)(x)
&\leq
\liminf_{j\to\infty}
(1-\alpha)I_\alpha(\mathcal D^\alpha w_j)(x)\\
&\leq
C_n\lim_{j\to\infty}
I_1(|\nabla w_j|)(x)\\
&=
C_n I_1(|\nabla w|)(x)\\
&<\infty. 
\end{aligned}
\]
Since \(E_w\) and the sequence \((w_j)\)
are independent of \(\alpha\), this
inequality holds for every \(x\in E_w\)
and every \(\alpha\in[1/2,1)\). This proves \eqref{claim-L5.1}. 

Finally, let
\(u,v\in W^{1,1}(\mathbb R^n)\) and set
\[
E_{u,v}=E_u\cap E_v\cap E_{u-v}.
\]
Since each of these sets has full measure, so does \(E_{u,v}\). Moreover, the reverse
triangle inequality yields, for every \(0<\alpha<1\),
\[
|\mathcal D^\alpha u(y)
-\mathcal D^\alpha v(y)|
\leq
\mathcal D^\alpha(u-v)(y)
\]
for almost every \(y\in\mathbb R^n\). Consequently, by \eqref{claim-L5.1}, for every \(x\in E_{u,v}\)
and every \(\alpha\in[1/2,1)\),
\[
\begin{aligned}
\bigl|
(1-\alpha)I_\alpha(\mathcal D^\alpha u)(x)
-
(1-\alpha)I_\alpha(\mathcal D^\alpha v)(x)
\bigr|
& \leq
(1-\alpha)I_\alpha(\mathcal D^\alpha(u-v))(x)\\
& \leq
C_n I_1(|\nabla(u-v)|)(x)\\
&<\infty.
\end{aligned}
\]
This proves \eqref{eq:uniform-stability-W11}.
\end{proof}

\subsection*{Proof of Theorem~\ref{thm:W11-full-limit}}

\begin{proof}

Let \(f\in W^{1,1}(\mathbb R^n)\). By the density of \(C_c^\infty(\mathbb R^n)\)
in \(W^{1,1}(\mathbb R^n)\), we can choose
\((f_j)\subset C_c^\infty(\mathbb R^n)\)
such that
\[
\|f_j-f\|_{W^{1,1}(\mathbb R^n)}
\leq 2^{-j},
\qquad j\geq1.
\]
Consequently,
\[
G=\sum_{j=1}^\infty
|\nabla(f_j-f)|
\in L^1(\mathbb R^n).
\]
As in the proof of Lemma~\ref{L5.1},
the set
\[
E_0=
\left\{
x\in\mathbb R^n:
I_1(G+|\nabla f|)(x)<\infty
\right\}
\]
has full measure, and
\[
\sum_{j=1}^{\infty}
I_1(|\nabla(f_j-f)|)(x)
=
I_1G(x)<\infty
\]
for every \(x\in E_0\).  Consequently,
\begin{equation}\label{eq:ae-gradient-approx}
I_1(|\nabla(f_j-f)|)(x)
\longrightarrow0
\qquad\text{as }j\to\infty
\end{equation}
for every \(x\in E_0\).

For each \(j\geq1\), let \(E_j\) be
the set of full measure provided by
Lemma~\ref{L5.1} for the pair \((f,f_j)\).
Set
\[
E=E_0\cap\bigcap_{j=1}^\infty E_j.
\]
Then \(E\) has full measure. Fix \(x\in E\). By Lemma~\ref{L5.1},
for every \(j\geq1\) and every
\(\alpha\in[1/2,1)\),
\[
\bigl|
(1-\alpha)I_\alpha(\mathcal D^\alpha f)(x)
-
(1-\alpha)I_\alpha(\mathcal D^\alpha f_j)(x)
\bigr|
\leq
C_n I_1(|\nabla(f-f_j)|)(x).
\]
Moreover,
\[
\bigl|
I_1(|\nabla f_j|)(x)
-
I_1(|\nabla f|)(x)
\bigr|
\leq
I_1(|\nabla(f_j-f)|)(x).
\]
Combining these inequalities, we obtain
\begin{equation}\label{final-ineq}
\begin{aligned}
\bigl|
(1-\alpha)I_\alpha(\mathcal D^\alpha f)(x)
-
K_n I_1(|\nabla f|)(x)
\bigr|
& \leq
\bigl|
(1-\alpha)I_\alpha(\mathcal D^\alpha f_j)(x)
-
K_n I_1(|\nabla f_j|)(x)
\bigr|\\
&\quad +
(C_n+K_n)
I_1(|\nabla(f_j-f)|)(x).
\end{aligned}
\end{equation}
For each fixed \(j\geq1\), since
\(f_j\in C_c^\infty(\mathbb R^n)\),
Theorem~\ref{thm:first-order-expansion} yields
\[
\lim_{\alpha\to1^-}
(1-\alpha)I_\alpha(\mathcal D^\alpha f_j)(x)
=
K_n I_1(|\nabla f_j|)(x).
\]
Hence, taking the limit superior as
\(\alpha\to1^-\) in \eqref{final-ineq},
we obtain
\[
\limsup_{\alpha\to1^-}
\bigl|
(1-\alpha)I_\alpha(\mathcal D^\alpha f)(x)
-
K_n I_1(|\nabla f|)(x)
\bigr|
\leq
(C_n+K_n)
I_1(|\nabla(f_j-f)|)(x).
\]
Finally, letting \(j\to\infty\)
and using \eqref{eq:ae-gradient-approx}, we conclude that
\[
\lim_{\alpha\to1^-}
(1-\alpha)I_\alpha(\mathcal D^\alpha f)(x)
=
K_n I_1(|\nabla f|)(x).
\]
Since \(E\) has full measure, the desired
almost-everywhere convergence follows.
\end{proof}

\section{Further remarks on the sharp constants}\label{sec6}

In this section, we discuss two further aspects
of the sharp-constant problem. We first examine
the difficulty of determining the unrestricted
optimal constant in dimension two for arbitrary
\(0<\alpha<1\). We then derive an explicit
convergent series representation for the sharp
constant \(C_{n,\alpha}^{*}\) in dimensions
\(n\geq3\).

\subsection{The unrestricted sharp constant in dimension two}

Although Theorems~\ref{thm:leading-asymptotics-2d}
and~\ref{thm:sharp-constant-radial-2d}
provide a precise description of the sharp-constant
problem in dimension two, the exact determination
of \(C_{2,\alpha}^{*}\) for arbitrary
\(0<\alpha<1\) remains unresolved by our approach.
Indeed, Lemmas~\ref{L2.1} and~\ref{L2.3} yield
\[
C_{2,\alpha}^{\mathrm{rad}}
=
(1-\alpha)
\frac{\gamma_{2,\alpha}}{\gamma_{2,1}}
A_\alpha(0)
\leq C_{2,\alpha}^{*}
\leq
(1-\alpha)
\frac{\gamma_{2,\alpha}}{\gamma_{2,1}}
\sup_{\beta\in[0,\pi/2]}A_\alpha(\beta).
\]
By Lemma~\ref{L2.2},
\[
A_\alpha(\beta)
=
\frac{4}{1-\alpha}+O(1)
\]
uniformly for \(\beta\in[0,\pi/2]\)
as \(\alpha\to1^-\).
Since
\[
\frac{\gamma_{2,\alpha}}{\gamma_{2,1}}
=
1+O(1-\alpha)
\qquad\text{as }\alpha\to1^-,
\]
we obtain
\[
0\leq
C_{2,\alpha}^{*}
-
C_{2,\alpha}^{\mathrm{rad}}
=
O(1-\alpha)
\qquad\text{as }\alpha\to1^-.
\]
Thus, the radial class is asymptotically
optimal at leading order, even though
its global optimality for each fixed
\(\alpha\in(0,1)\) has not been established.

A natural question is whether
\[
C_{2,\alpha}^{*}
=
C_{2,\alpha}^{\mathrm{rad}},
\qquad 0<\alpha<1.
\]
A sufficient condition would be
\[
A_\alpha(\beta)\leq A_\alpha(0),
\qquad \beta\in[0,\pi/2],
\]
which would identify the radial direction
as a maximizer of the associated directional
functional. However, the monotonicity argument
used in the proof of Lemma~\ref{L3.4}
for \(n\geq3\) does not extend directly
to dimension two. Indeed, that argument
relies on
\[
a=\frac{n-\alpha}{2}>1
\]
and on the additional term
\[
-\frac{n-2}{3}\tau
\]
in the estimate for the logarithmic derivative
of the transformed kernel. In dimension two,
we instead have
\[
a=\frac{2-\alpha}{2}<1,
\qquad n-2=0,
\]
so the same reasoning no longer establishes
the required monotonicity.

It remains an interesting question whether
\[
C_{2,\alpha}^{*}
=
C_{2,\alpha}^{\mathrm{rad}}
\]
for every \(0<\alpha<1\), or at least for
\(\alpha\) sufficiently close to \(1\).

\subsection{An explicit series representation in higher dimensions}

We now complement the integral representation
\eqref{eq:def-radial-coefficient-n} with an
absolutely convergent series for
\(\mathfrak C_{n,\alpha}^{\rm rad}\).
Combined with \eqref{eq:sharp-Cn}, this yields
an explicit series representation for the
sharp constant \(C_{n,\alpha}^{*}\) for every
\(n\geq3\) and \(0<\alpha<1\).

\begin{proposition}\label{prop:series-Cn}
Let \(n\geq3\) and \(0<\alpha<1\). Then
\begin{equation}\label{eq:series-radial-Cn}
\mathfrak C_{n,\alpha}^{\rm rad}
=
|\mathbb S^{n-1}|
\sum_{k=0}^\infty
\frac{
\left(\frac{n+\alpha}{2}\right)_k
\left(1+\frac{\alpha}{2}\right)_k
}{
\left(\frac n2\right)_k\,k!
}
\left(
\frac{1}{(\alpha+2k)^2}
+
\frac{1}{(n+2k)^2}
\right).
\end{equation}
The series converges absolutely for every
\(0<\alpha<1\).
\end{proposition}

\begin{proof}

Fix \(0<\alpha<1\) and set
\[
p=\frac{n+\alpha}{2}.
\]
By \eqref{eq:def-K-alpha-n}, for
\(\tau>0\) and \(0\leq\theta\leq\pi\),
we have
\[
\mathscr K_{n,\alpha}(\tau,\theta)
=
\frac{e^{-\alpha\tau}+e^{-n\tau}}
{(1-2e^{-\tau}\cos\theta+e^{-2\tau})^p}.
\]
By \eqref{eq:def-radial-coefficient-n}, 
\[
\mathfrak C_{n,\alpha}^{\rm rad}
=
|\mathbb S^{n-2}|
\int_0^\infty
\tau(e^{-\alpha\tau}+e^{-n\tau})\int_0^\pi
\frac{(\sin\theta)^{n-2}}
{(1-2e^{-\tau}\cos\theta+e^{-2\tau})^p}
\,d\theta\,d\tau.
\]
Making the change of variables \(t=e^{-\tau}\),
we obtain
\[
\mathfrak C_{n,\alpha}^{\rm rad}
=
|\mathbb S^{n-2}|
\int_0^1
(-\log t)(t^{\alpha-1}+t^{n-1})
\int_0^\pi
\frac{(\sin\theta)^{n-2}}
{(1-2t\cos\theta+t^2)^p}
\,d\theta\,dt.
\]
Now, using the classical hypergeometric identity
(see, e.g., \cite[Sections~2.2 and~3.1]{AndrewsAskeyRoy})
\[
|\mathbb S^{n-2}|
\int_0^\pi
\frac{(\sin\theta)^{n-2}}
{(1-2t\cos\theta+t^2)^p}
\,d\theta
=
|\mathbb S^{n-1}|
{}_2F_1\left(
p,1+\frac{\alpha}{2};
\frac n2;t^2
\right),
\qquad 0<t<1,
\]
we deduce that
\[
\mathfrak C_{n,\alpha}^{\rm rad}
=
|\mathbb S^{n-1}|
\int_0^1
(-\log t)(t^{\alpha-1}+t^{n-1})
{}_2F_1\left(
p,1+\frac{\alpha}{2};
\frac n2;t^2
\right)\,dt.
\]
Expanding the hypergeometric function
into its defining series
(see, e.g., \cite[Section~2.1]{AndrewsAskeyRoy})
and using
\[
\int_0^1(-\log t)t^{\lambda-1}\,dt
=
\frac1{\lambda^2},
\qquad \lambda>0,
\]
we obtain \eqref{eq:series-radial-Cn},
where the interchange of summation and
integration is justified by the
nonnegativity of all terms.

Finally, by the standard asymptotic properties
of the Gamma function, we have
\[
\frac{
\left(\frac{n+\alpha}{2}\right)_k
\left(1+\frac{\alpha}{2}\right)_k
}{
\left(\frac n2\right)_k\,k!
}
=
O(k^\alpha)
\qquad\text{as }k\to\infty.
\]
Thus, the summands in
\eqref{eq:series-radial-Cn} are
\(O(k^{\alpha-2})\) as \(k\to\infty\).
Since \(0<\alpha<1\), the series
converges absolutely.
\end{proof}

By Legendre's duplication formula,
\[
\Gamma\left(\frac{n-1}{2}\right)
\Gamma\left(\frac n2\right)
=
2^{2-n}\sqrt{\pi}\,\Gamma(n-1),
\]
and \eqref{eq:gamma-n-beta}, we obtain
\[
\frac{\gamma_{n,\alpha}}{\gamma_{n,1}}
|\mathbb S^{n-1}|
=
\frac{
2^{n-\alpha}\pi^{n/2}
\Gamma\left(\frac{n-\alpha}{2}\right)
}{
\Gamma(\alpha/2)\Gamma(n-1)
}.
\]
Combining this identity with
\eqref{eq:sharp-Cn} and
Proposition~\ref{prop:series-Cn},
we obtain the following explicit series representation
of the sharp constant \(C_{n,\alpha}^{*}\).

\begin{corollary}
Let \(n\geq3\) and \(0<\alpha<1\). Then
\[
C_{n,\alpha}^{*}
=
\frac{
2^{n-\alpha}\pi^{n/2}(1-\alpha)
\Gamma\left(\frac{n-\alpha}{2}\right)
}{
\Gamma(\alpha/2)\Gamma(n-1)
}\sum_{k=0}^{\infty}
\frac{
\left(\frac{n+\alpha}{2}\right)_k
\left(1+\frac{\alpha}{2}\right)_k
}{
\left(\frac n2\right)_k\,k!
}
\left(
\frac{1}{(\alpha+2k)^2}
+
\frac{1}{(n+2k)^2}
\right).
\]
The series converges absolutely.
\end{corollary}

\vspace{1cm}

\noindent {\bf Data availability.} No data was used for the research described in the article.

\medskip

\noindent {\bf Conflict of interest.} The authors declare that they have no conflict of interest.

\end{document}